\documentclass[11pt]{amsart}

\usepackage[a4paper,margin=1.15in]{geometry}
\usepackage{microtype}

\usepackage{amsmath,amssymb,amsfonts,amsthm}
\usepackage{mathtools}
\usepackage{mathrsfs}
\usepackage{bm}

\usepackage{graphicx}
\usepackage{booktabs}
\usepackage{enumitem}

\usepackage[dvipsnames]{xcolor}
\usepackage[
    colorlinks=true,
    linkcolor=MidnightBlue,
    citecolor=ForestGreen,
    urlcolor=BrickRed
]{hyperref}
\usepackage[capitalize,nameinlink]{cleveref}

\numberwithin{equation}{section}

\theoremstyle{plain}
\newtheorem{theorem}{Theorem}[section]
\newtheorem{lemma}[theorem]{Lemma}

\newtheorem{corollary}[theorem]{Corollary}

\theoremstyle{definition}
\newtheorem{definition}[theorem]{Definition}

\theoremstyle{remark}
\newtheorem{remark}[theorem]{Remark}

\crefname{theorem}{Theorem}{Theorems}
\crefname{lemma}{Lemma}{Lemmas}
\crefname{proposition}{Proposition}{Propositions}
\crefname{corollary}{Corollary}{Corollaries}
\crefname{definition}{Definition}{Definitions}
\crefname{assumption}{Assumption}{Assumptions}
\crefname{remark}{Remark}{Remarks}
\crefname{equation}{equation}{equations}
\crefname{section}{Section}{Sections}

\usepackage{algorithm}    % 算法排版环境
\usepackage{algorithmic}  % 算法伪代码排版（旧版）

\title[LbNM for biharmonic problems]{A Learning-Based Numerical Method for Biharmonic boundary value problems via  the Almansi Decomposition}

\author{Mengjia Bai}
\address{School of Mathematical Sciences, Fudan University}
\email{bai\_mj@fudan.edu.cn}

\begin{document}

\begin{abstract}

We extend the learning-based numerical method (LbNM) to two-dimensional biharmonic boundary value problems, combining high accuracy with efficient repeated solves for different boundary data on a fixed domain. The method reconstructs a discrete boundary-to-solution operator from explicit biharmonic samples generated by the Almansi decomposition, using their boundary values and normal derivatives through Tikhonov regularization.
Once constructed at prescribed evaluation points, the operator provides solutions for new boundary inputs by a single matrix-vector multiplication. Moreover, additional sample solutions can be incorporated through rank-one updates without rebuilding the operator. We establish a quantitative convergence estimate under a biharmonic continuation assumption. Numerical experiments on nonconvex domains demonstrate accuracy close to the double-precision limit for a smooth benchmark and high accuracy for solutions induced by nearby exterior singularities.

% In this paper, we extend the learning-based numerical method (LbNM) to two-dimensional biharmonic boundary value problems with prescribed boundary values and normal derivatives. Using the Almansi decomposition, we construct explicit biharmonic samples from harmonic polynomial modes and their radial multiples. Their boundary traces are used to reconstruct a discrete boundary-to-solution operator through Tikhonov regularization. Since the samples satisfy the governing equation exactly, no interior discretization is required. Once constructed at prescribed evaluation points, the operator provides solutions for different boundary inputs through matrix-vector multiplication. With the boundary discretization and regularization parameter fixed, it can also be updated online as additional sample solutions become available. We establish a quantitative error estimate under a biharmonic continuation assumption. Numerical experiments on nonconvex domains demonstrate high accuracy and low computational cost for repeated solution evaluation.


\end{abstract}

\maketitle
%\tableofcontents %目录

\section{Introduction}
%问题背景 → 数学观察 → 最相关已有方法 → LbNM 发展 → 本文算法 → 方法特点与贡献 → 文章结构

The biharmonic equation is a classical fourth-order elliptic  equation arising in thin-plate bending theory \cite{Timoshenko_plates}, plane elasticity \cite{Timoshenko_elasticity}, and slow viscous flows \cite{Happel_stokes}. In this work, we consider the boundary value problem 
\begin{equation*}
\Delta^2 u = 0 \quad \text{in } \Omega,
\qquad
u = f,\quad \partial_{n}u = g
\quad \text{on } \partial\Omega,
\end{equation*}
where $\Omega\subset\mathbb{R}^2$ is a bounded domain and $n$ denotes the outward unit normal vector on $\partial\Omega$.
In standard conforming finite element formulations, the requirement of $H^2$ conformity complicates the construction of approximation spaces, while the simultaneous treatment of the two boundary conditions adds to the implementation effort \cite{engel2002continuous,brenner2005c}. When different boundary profiles are considered on the same domain, it is desirable to combine high accuracy with a low computational cost for repeated solves.

Under suitable regularity assumptions on the domain and the boundary data, the solution satisfies an estimate of the form~\cite{Gazzola_polyharmonic,grisvard2011elliptic}
\begin{equation*}
\|u\|_{H^2(\Omega)}
\leq
C\left(
\|f\|_{H^{3/2}(\partial\Omega)}
+
\|g\|_{H^{1/2}(\partial\Omega)}
\right).
\end{equation*}
Based on the stated well-posedness and regularity assumptions, the interior solution is uniquely determined by the prescribed boundary data. Equivalently, the problem can be viewed in terms of a boundary-to-solution operator
\begin{equation*}
\mathcal{S}:
H^{3/2}(\partial\Omega)\times H^{1/2}(\partial\Omega)
\longrightarrow H^2(\Omega),
\qquad
\mathcal{S}(f,g)=u.
\end{equation*}
This provides a natural starting point for methods based on boundary information. The method of fundamental solutions (MFS), for example, constructs approximations from functions that satisfy the biharmonic equation exactly in the physical domain \cite{karageorghis1987method,li2011trefftz,dou2018analysis}. By placing their singularities outside the domain, the MFS avoids interior discretization and reduces the computation to enforcing the boundary conditions. It can achieve high accuracy for smooth solutions, but its performance may depend on the placement of source points and the conditioning of the resulting system \cite{dou2018analysis,gaspar2025method,antunes2026improving}. Alternative formulations and stabilization strategies have been developed to address these issues \cite{gaspar2025method,antunes2026improving}.

The Learning-Based Numerical Method (LbNM), introduced by Chen et al. \cite{chen2025learning} for high-frequency Helmholtz equations, constructs a regularized discrete solution operator from representative sample solutions. Related developments include acoustic frequency-domain simulations with mixed boundary conditions \cite{li2024acoustic} and a formulation based on Bessel functions for dissipative media and complex geometries \cite{song2026robust}. Extending this framework to the biharmonic equation requires a suitable sample space and an analysis of the two prescribed boundary traces. The Almansi decomposition suggests constructing samples from two harmonic families. However, both families contribute to the boundary values and normal derivatives, so their approximation properties must be connected to the simultaneous fitting of these traces and to the resulting interior error.

In this paper, we extend LbNM for the two-dimensional biharmonic boundary value problems, using the Almansi decomposition
$u = \phi + r^2\psi$, 
with $\phi$ and $\psi$ are harmonic functions. We construct the samples from harmonic polynomial modes, normalized on the boundary of an enclosing fictitious disk, and their multiples by $r^2$. These explicit functions satisfy the biharmonic equation without introducing artificial singular source points. Their boundary values and normal derivatives are assembled into a sample matrix, from which the discrete solution operator is reconstructed through Tikhonov regularization. The reconstruction admits an equivalent formulation as a regularized least-squares fit in the biharmonic sample space. Since the governing equation is satisfied by construction, no interior discretization or evaluation of a fourth-order PDE residual is required.

The construction retains the operator reuse and incremental update mechanisms of LbNM. For a fixed domain, discretization, and regularization parameter, the reconstructed operator can be applied to different boundary inputs without being rebuilt. Once it is stored at the prescribed evaluation points, each subsequent solve reduces to matrix-vector multiplication. The sample collection can also be enriched online: with the boundary collocation set fixed, adding a new sample solution gives a rank-one update of the sample cross-product matrix. Thus, the new sample information can be incorporated without reassembling the contributions of preceding samples.

We establish a quantitative error estimate for the resulting biharmonic reconstruction. Under a biharmonic continuation assumption on a larger concentric disk and a suitable choice of the regularization parameter, the $L^2(\Omega)$  error is bounded by a geometrically decaying truncation term and a boundary collocation term of order $N^{-1}$ , where $N$ is the number of boundary collocation points. The bound is expressed directly in terms of the $H^2$ norm of the extended exact solution.
Numerical experiments on nonconvex geometries, including tests with solutions induced by exterior singularities located close to the computational boundary, demonstrate the accuracy and robustness of the proposed method.

The remainder of this paper is organized as follows. Section \ref{sec:learning_based_numerical_method} introduces the boundary-to-solution formulation, presents the LbNM reconstruction algorithm, and describes the construction of the sample solutions. Section \ref{sec:theory} develops the approximation theory for the harmonic bases and the Almansi decomposition and establishes the convergence estimate for the proposed method. Section \ref{sec:experiments} reports numerical experiments for smooth solutions, complex geometries, and solutions with nearby  singularities. Finally, Section \ref{sec:conclusions} summarizes the main results and discusses possible extensions.

\section{The Learning-Based Numerical Method}
\label{sec:learning_based_numerical_method}
We now introduce the learning-based numerical method for the two-dimensional biharmonic equation. Let $\Omega\subset\mathbb{R}^2$ be a bounded domain with a sufficiently smooth boundary. We consider the boundary value problem:
\begin{equation}\label{eq:biharmonic_cases}
\begin{cases}
    \Delta^2 u = 0, & \text{in } \Omega, \\
    u = f,\  \partial_{n} u = g, & \text{on } \partial\Omega,
\end{cases}
\end{equation}
where $n$  denotes the outward unit normal vector. The proposed method reconstructs a finite-dimensional
approximation of the boundary-to-solution operator from a collection of biharmonic sample solutions. We first motivate this operator viewpoint, then present the regularized reconstruction procedure, and finally describe the construction of the sample space through the Almansi decomposition.

\subsection{Boundary-to-solution formulation}
Suppose that $u$ solves \eqref{eq:biharmonic_cases}. For a fixed point $\boldsymbol{x} \in \Omega$, let $G(\boldsymbol{x},\boldsymbol{y})$ denote the clamped biharmonic Green function, with $\boldsymbol{y}$ as the integration variable, satisfying
\begin{equation}\label{eq:green_function}
\begin{cases}
    \Delta^2 G(\boldsymbol{x}, \boldsymbol{y}) = \delta(\boldsymbol{y} - \boldsymbol{x}), & \boldsymbol{y} \in \Omega, \\
    G(\boldsymbol{x}, \boldsymbol{y}) = 0, \  \partial_{n} G(\boldsymbol{x}, \boldsymbol{y}) = 0, & \boldsymbol{y} \in \partial\Omega,
\end{cases}
\end{equation}
where $\delta(\cdot)$ denotes the Dirac delta distribution concentrated at $\boldsymbol{x}$. Applying Green's identity for the biharmonic operator and using the homogeneous boundary conditions for $G$ gives
\begin{equation*}
\begin{aligned}
    \int_{\Omega} \big( (\Delta^2 u)G - u(\Delta^2 G) \big) \, \mathrm{d}\boldsymbol{y}
    &= \int_{\partial\Omega} \big( \partial_{n}(\Delta u)G - \Delta u\partial_{n}G + \Delta G \partial_{n}u - \partial_{n}(\Delta G)u \big) \, \mathrm{d}S_{\boldsymbol{y}} \\
    &= \int_{\partial\Omega} \big( \Delta G \partial_{n}u - \partial_{n}(\Delta G)u \big) \, \mathrm{d}S_{\boldsymbol{y}}.
\end{aligned}
\end{equation*}
Substituting the prescribed boundary data into this identity yields the representation
\begin{equation}
\begin{aligned}
\label{eq:green_representation}
    u(\boldsymbol{x}) 
    &= \int_{\partial\Omega} \big( \partial_{n}(\Delta G(\boldsymbol{x}, \boldsymbol{y})) u(\boldsymbol{y}) - \Delta G(\boldsymbol{x}, \boldsymbol{y}) \partial_{n} u(\boldsymbol{y}) \big) \, \mathrm{d}S_{\boldsymbol{y}} \\
    &= \int_{\partial\Omega} \big( \partial_{n}(\Delta G(\boldsymbol{x}, \boldsymbol{y})) f(\boldsymbol{y}) - \Delta G(\boldsymbol{x}, \boldsymbol{y}) g(\boldsymbol{y}) \big) \, \mathrm{d}S_{\boldsymbol{y}}.
\end{aligned}
\end{equation}
Equation~\eqref{eq:green_representation} shows that the interior solution is
determined by the prescribed boundary data. We therefore introduce the
boundary-to-solution operator
\begin{equation}
\label{eq:solution_operator}
\mathcal{A}:
H^{3/2}(\partial\Omega)\times H^{1/2}(\partial\Omega)
\longrightarrow H^2(\Omega),
\qquad
\mathcal{A}(f,g)=u.
\end{equation}
Formally, this operator may be decomposed as
\begin{equation}
\label{eq:operator_decomposition}
\mathcal{A}(f,g)(\boldsymbol{x})
=
\mathcal{A}^{(1)}f(\boldsymbol{x})
+
\mathcal{A}^{(2)}g(\boldsymbol{x}),
\end{equation}
where
\begin{equation*}
\mathcal{A}^{(1)}f(\boldsymbol{x})
=
\int_{\partial\Omega}
\partial_{n}
\bigl(\Delta G(\boldsymbol{x},\boldsymbol{y})\bigr)f(\boldsymbol{y})\,\mathrm d S_{\boldsymbol{y}}, \quad
\mathcal{A}^{(2)}g(\boldsymbol{x})=
-\int_{\partial\Omega}
\Delta G(\boldsymbol{x},\boldsymbol{y})g(\boldsymbol{y})\,\mathrm d S_{\boldsymbol{y}}.
\end{equation*}

For a finite-dimensional discretization, suppose that the boundary data are approximated by
\begin{equation*}
    f \approx \sum_{i=1}^N f(\boldsymbol{x}_i) e_{i}(\boldsymbol{x}), \quad
    g \approx \sum_{i=1}^N g(\boldsymbol{x}_i) e_{i}(\boldsymbol{x}),
\end{equation*}
where $\{e_{i}\}_{i=1}^N$ are boundary interpolation functions. The corresponding approximation  $u^{(c)}$ then takes the form
\begin{equation*}
\begin{aligned}
    u^{(c)}(\boldsymbol{x})
    &= \mathcal{A}^{(1)}\Big(\sum_{i=1}^N f(\boldsymbol{x}_i) e_i(\boldsymbol{x}) \Big) + \mathcal{A}^{(2)}\Big(\sum_{i=1}^N g(\boldsymbol{x}_i) e_i(\boldsymbol{x}) \Big) \\
    &= \sum_{i=1}^{N} \big( f(\boldsymbol{x}_i) \mathcal{A}^{(1)}e_i(\boldsymbol{x}) + g(\boldsymbol{x}_i) \mathcal{A}^{(2)}e_i(\boldsymbol{x}) \big) 
    = \sum_{i=1}^{N} \begin{pmatrix} \mathcal{A}^{(1)}_i(\boldsymbol{x}), \mathcal{A}^{(2)}_i(\boldsymbol{x}) \end{pmatrix} \begin{pmatrix} f(\boldsymbol{x}_i) \\ g(\boldsymbol{x}_i) \end{pmatrix},
\end{aligned}
\end{equation*}
where $\mathcal{A}_i^{(k)}(\boldsymbol{x}) = \mathcal{A}^{(k)} e_i(\boldsymbol{x})$ for $k=1,2$.
The functions $\{\mathcal{A}^{(1)}_j,\mathcal{A}^{(2)}_j\}_{j=1}^N$ constitute a discrete representation of the boundary-to-solution operator. They are generally unavailable in closed form because the Green function depends on the geometry of $\Omega$. The central idea of the proposed LbNM is therefore to reconstruct these functions from a collection of  biharmonic sample solutions, without explicitly computing the Green function.

\subsection{Algorithm}
We next describe the regularized reconstruction of the discrete solution operator from sample observations. 

Let $\{u_i\}_{i=1}^M$ be biharmonic sample solutions associated with boundary-data pairs $(f_i, g_i)$ and let $T=\{\boldsymbol{x}_j\}_{j=1}^N \subset \partial\Omega$ be the boundary collocation set. For each sample solution and each evaluation point $\boldsymbol{x} \in \Omega$, we seek an  approximation of the form
\begin{equation*}
\begin{aligned}
    u_i(\boldsymbol{x}) 
    \approx \sum_{j=1}^N \begin{pmatrix} A^{(1)}_j(\boldsymbol{x}), A^{(2)}_j(\boldsymbol{x}) \end{pmatrix} \begin{pmatrix} f_{i,j} \\ g_{i,j} \end{pmatrix} 
    = \sum_{j=1}^N \begin{pmatrix} A^{(1)}_j(\boldsymbol{x}), A^{(2)}_j(\boldsymbol{x}) \end{pmatrix} \begin{pmatrix} u_i(\boldsymbol{x}_j) \\ \partial_{n}u_i(\boldsymbol{x}_j) \end{pmatrix}.
\end{aligned}
\end{equation*}

Consequently, assembling the approximations for all $M$ training samples yields the following linear system:
\begin{equation*}
\begin{aligned}
    &\begin{pmatrix} u_1(\boldsymbol{x}) & \cdots & u_M(\boldsymbol{x}) \end{pmatrix} \\
    &= \begin{pmatrix} A_1^{(1)}(\boldsymbol{x}) & A_1^{(2)}(\boldsymbol{x}) & \cdots & A_N^{(1)}(\boldsymbol{x}) & A_N^{(2)}(\boldsymbol{x}) \end{pmatrix}
    \begin{pmatrix} 
        u_1(\boldsymbol{x}_1) & \cdots & u_M(\boldsymbol{x}_1) \\
        \partial_{n} u_1(\boldsymbol{x}_1) & \cdots & \partial_{n} u_M(\boldsymbol{x}_1) \\
        \vdots & \ddots & \vdots \\
        u_1(\boldsymbol{x}_N) & \cdots & u_M(\boldsymbol{x}_N) \\
        \partial_{n} u_1(\boldsymbol{x}_N) & \cdots & \partial_{n} u_M(\boldsymbol{x}_N)
    \end{pmatrix},
\end{aligned}
\end{equation*}
which can be compactly denoted in matrix form as
\begin{equation*}
    \boldsymbol{b}(\boldsymbol{x})_{1\times M} = \boldsymbol{a}(\boldsymbol{x})_{1\times 2N} \boldsymbol{V}_{2N \times M}.
\end{equation*}
Depending on the relative sizes of $M$ and $2N$, this system may be overdetermined, square, or underdetermined. The regularized formulation below is well defined in all three cases. We obtain the optimal basis operator vector $\boldsymbol{a}_*$ by solving the following Tikhonov minimization problem:
\begin{equation}\label{eq:tikhonov_row}
    \boldsymbol{a}_* = \mathop{\operatorname{argmin}}_{\boldsymbol{a} \in \mathbb{C}^{1 \times 2N}} \left\{ \|\boldsymbol{a}\boldsymbol{V} - \boldsymbol{b}\|_2^2 + \alpha \|\boldsymbol{a}\|_2^2 \right\},
\end{equation}
where $\alpha \|\boldsymbol{a}\|_2^2$ is the penalty term with an a priori given regularization parameter $\alpha > 0$.
Then, there exists a unique minimizer $\boldsymbol{a}_*$ for any $\alpha > 0$, which is explicitly given by
\begin{equation}\label{eq:solution_operator}
    \boldsymbol{a}_* = \boldsymbol{b}\boldsymbol{V}^*(\boldsymbol{V}\boldsymbol{V}^* + \alpha \boldsymbol{I})^{-1},
\end{equation}
Here, $\boldsymbol{I}$ denotes the $2N \times 2N$ identity matrix, and $\boldsymbol{V}^*$ represents the conjugate transpose of $\boldsymbol{V}$. 

For any new boundary input pair $(f, g)$, define the complex reconstruction 
\begin{equation*}
    \widetilde{u}^{(c)}(\boldsymbol{x}) \approx \boldsymbol{a}_* \boldsymbol{h} = \boldsymbol{b}\boldsymbol{V}^*(\boldsymbol{V}\boldsymbol{V}^* + \alpha \boldsymbol{I})^{-1}\boldsymbol{h},
\end{equation*}
where $\boldsymbol{h} = \big(f(\boldsymbol{x}_1), g(\boldsymbol{x}_1), \cdots, f(\boldsymbol{x}_N), g(\boldsymbol{x}_N)\big)^T \in \mathbb{R}^{2N \times 1}$ is the vector of nodal boundary data. For real-valued problems, $\boldsymbol{h}$ is real, while the complex notation
remains convenient for the harmonic modes. The reported real-valued approximation is $u^{(c)} = \operatorname{Re}\widetilde u^{(c)} $

To describe the incorporation of additional samples, we make the dependence on the sample size explicit in this paragraph by writing 
\begin{equation*}
    \boldsymbol{V}_M := \boldsymbol{V}, \quad
    \boldsymbol{b}_M(\boldsymbol{x})
    := \boldsymbol{b}(\boldsymbol{x}), \quad
    \boldsymbol{a}_{*,M}(\boldsymbol{x})
    := \boldsymbol{a}_{*}(\boldsymbol{x}).
\end{equation*}
The boundary collocation set $T$, the existing sample, and the regularization parameter $\alpha>0$ are kept fixed during the update. Define
\begin{equation}\label{eq:sigma_definition}
    \boldsymbol{\Sigma}_M := \boldsymbol{V}_M\boldsymbol{V}^*_M,\quad
    \boldsymbol{P}_M 
    =
    (\boldsymbol{\Sigma}_M + \alpha\boldsymbol{I})^{-1}.
\end{equation}
The reconstructed operator is then
\begin{equation*}
  \boldsymbol{a}_{*,M}(\boldsymbol{x})
  =
  \boldsymbol{b}_{M}(\boldsymbol{x})\boldsymbol{V}^*_M\boldsymbol{P}_M. 
\end{equation*}
The $(i,j)$ entry of $\boldsymbol{\Sigma}_M$ can be written explicitly as
\begin{equation*}
(\boldsymbol{\Sigma}_M)_{i,j} =
\begin{cases}
    \sum_{\ell=1}^M u_\ell(\boldsymbol{x}_m) \overline{u_\ell(\boldsymbol{x}_n)}, & \text{if } i = 2m-1, \ j = 2n-1, \\[1.5ex]
    \sum_{\ell=1}^M \partial_{n} u_\ell(\boldsymbol{x}_m) \overline{u_\ell(\boldsymbol{x}_n)}, & \text{if } i = 2m, \qquad j = 2n-1, \\[1.5ex]
    \sum_{\ell=1}^M u_\ell(\boldsymbol{x}_m) \overline{\partial_{n} u_\ell(\boldsymbol{x}_n)}, & \text{if } i = 2m-1, \ j = 2n, \\[1.5ex]
    \sum_{\ell=1}^M \partial_{n} u_\ell(\boldsymbol{x}_m) \overline{\partial_{n} u_\ell(\boldsymbol{x}_n)}, & \text{if } i = 2m, \qquad j = 2n,
\end{cases}
\end{equation*}
where $1 \le m, n \le N$.

Now, upon introducing a new sample solution $u_{M+1}(\boldsymbol{x})$, we capture its boundary evaluations using the column vector $\boldsymbol{v}_{M+1} \in \mathbb{C}^{2N \times 1}$:
\begin{equation*}
    \boldsymbol{v}_{M+1} = \begin{pmatrix} u_{M+1}(\boldsymbol{x}_1) & \partial_{n} u_{M+1}(\boldsymbol{x}_1) & \cdots & u_{M+1}(\boldsymbol{x}_N) & \partial_{n} u_{M+1}(\boldsymbol{x}_N) \end{pmatrix}^T.
\end{equation*}
Consequently, denoting the updated matrix by $\boldsymbol{\Sigma}_{M+1} := \boldsymbol{V}_1\boldsymbol{V}_1^*$, we obtain the rank-one update
\begin{equation}\label{eq:rank_one_update}
    \boldsymbol{\Sigma}_{M+1} = \boldsymbol{\Sigma}_M + \boldsymbol{v}_{M+1} \boldsymbol{v}_{M+1}^*.
\end{equation}
Thus, when samples arrive sequentially or are generated adaptively, \eqref{eq:rank_one_update} incorporates the new boundary trace through one outer product and avoids recomputing the contributions of the preceding samples.

Since $\alpha>0$, the matrix $P_M$ is Hermitian positive definite,
and hence
\begin{equation*}
    1+\boldsymbol{v}_{M+1}^*P_M\boldsymbol{v}_{M+1}>0.
\end{equation*}
Direct multiplication gives the inverse update
\begin{equation}
\label{eq:inverse-rank-one-update}
    \boldsymbol{P}_{M+1}
    =
    \boldsymbol{P}_M
    -
    \frac{
        \boldsymbol{P}_M\boldsymbol{v}_{M+1}
        \boldsymbol{v}_{M+1}^*\boldsymbol{P}_M
    }{
        1+\boldsymbol{v}_{M+1}^*\boldsymbol{P}_M\boldsymbol{v}_{M+1}
    }.
\end{equation}
Moreover,
\begin{equation*}
    \boldsymbol{b}_{M+1}(x)\boldsymbol{V}_{M+1}^*
    =
    \boldsymbol{b}_M(x)\boldsymbol{V}_M^*
    +
    u_{M+1}(x)\boldsymbol{v}_{M+1}^*.
\end{equation*}
Substituting these identities into
\begin{equation*}
    \boldsymbol{a}_{*,M+1}(x)
    =
    \boldsymbol{b}_{M+1}(x)\boldsymbol{V}_{M+1}^*\boldsymbol{P}_{M+1}
\end{equation*}
yields
\begin{equation}
\label{eq:operator-rank-one-update}
\begin{aligned}
    \boldsymbol{a}_{*,M+1}(x)
    =
    \boldsymbol{a}_{*,M}(x)
    +
    \frac{
        \bigl(
            u_{M+1}(x)
            -
            \boldsymbol{a}_{*,M}(x)\boldsymbol{v}_{M+1}
        \bigr)
        \boldsymbol{v}_{M+1}^*\boldsymbol{P}_M
    }{
        1+\boldsymbol{v}_{M+1}^*\boldsymbol{P}_M\boldsymbol{v}_{M+1}
    }.
\end{aligned}
\end{equation}
Thus, both the inverse matrix and the discrete operator can be
updated using the new sample and the previously computed quantities, without reassembling the contributions of the existing samples.

\subsection{Construction of the  Sample Space}
We now specify the samples used in the reconstruction.

For a systematic scaling of the harmonic modes, we embed the physical domain $\Omega$ into a larger fictitious circular disk $O_\rho = \{ (r,\theta) : r < \rho \}$ such that $\overline{\Omega} \subset O_\rho$, where $\partial O_\rho$ denotes the boundary of the disk with radius $\rho$. Motivated by the classical Almansi decomposition, we construct the sample space from functions of the form
\begin{equation}\label{eq:almansi_new}
    u = \phi + r^2 \psi,
\end{equation}
where $\phi$ and $\psi$ are harmonic. The precise representation and its stability properties are established in Section~\ref{sec:theory}.

To approximate these harmonic components, we introduce a set of normalized harmonic polynomial basis functions defined on the boundary $\partial O_\rho$:
\begin{equation}\label{eq:normalized_basis}
    \varphi_n(r,\theta) = \frac{r^{|n|}e^{\mathrm{i}n\theta}}{\|r^{|n|}e^{\mathrm{i}n\theta}\|_{L^2(\partial O_{\rho})}} = \frac{r^{|n|}e^{\mathrm{i}n\theta}}{\sqrt{2\pi\rho}\rho^{|n|}}, \quad n = -K, \cdots, K,
\end{equation}
where $\mathrm{i}$ is the imaginary unit, and $K$ is the maximum truncation degree. Using \eqref{eq:almansi_new}, we construct the finite-dimensional analytical sample family from the two sets
\begin{equation*}
    \mathcal{B}^{(1)} = \big\{ \varphi_n(r,\theta) \big\}_{n=-K}^K \quad \text{and} \quad \mathcal{B}^{(2)} = \big\{ r^2 \varphi_n(r,\theta) \big\}_{n=-K}^K.
\end{equation*}
Consequently, the construction yields a total of $M = 4K+2$  analytical training solutions. Because the modes with indices $n$ and $-n$ occur in complex-conjugate pairs, real-valued approximations can be recovered by taking the real part of the reconstructed solution.

\begin{algorithm}[H]
	\caption{Learning-Based Numerical Method (LbNM) with Tikhonov Regularization}
	\label{alg:LbNM_Biharmonic}
	\begin{algorithmic}[1]
		\STATE \textbf{Input:} 
        Computational domain $\Omega$;
        Boundary collocation points $T = \{\boldsymbol{x}_j\}_{j=1}^N \subset \partial\Omega$; 
		truncation degree $K$; 
        fictitious-disk radius $\rho$ with $\overline{\Omega} \subset O_\rho$;
		regularization parameter $\alpha > 0$; 
		evaluation point $\boldsymbol{x} \in \Omega$.
		\STATE \textbf{Output:} Numerical solution $u^{(c)}(\boldsymbol{x})$.
		
		\STATE \COMMENT{Step 1: Basis Functions Construction (Training Solutions)}
		\STATE Initialize $\mathcal{U} = \emptyset$.
		\FOR{$n = - K$ \TO $K$}
			\STATE Generate harmonic mode: $u_{n,1}(r, \theta) = \varphi_n(r,\theta)$
			\STATE Generate biharmonic mode: $u_{n,2}(r, \theta) = r^{2} \varphi_n(r,\theta)$
			\STATE $\mathcal{U} \leftarrow \mathcal{U} \cup \{u_{n,1}, u_{n,2}\}$
		\ENDFOR
		\STATE Total samples $M = 4K + 2$.
		
		\STATE \COMMENT{Step 2: Operator Matrix Assembly}
		\STATE Initialize $\boldsymbol{V} \in \mathbb{C}^{2N \times M}$.
		\FOR{$\ell = 1$ \TO $M$}
			\STATE Evaluate $u_\ell$ and $\partial_{n} u_\ell$ at $T$.
			\STATE Set $\ell$-th column $\boldsymbol{v}_\ell = [u_\ell(\boldsymbol{x}_1), \partial_{n} u_\ell(\boldsymbol{x}_1), \dots, u_\ell(\boldsymbol{x}_N), \partial_{n} u_\ell(\boldsymbol{x}_N)]^T$.
		\ENDFOR
		
		\STATE \COMMENT{Step 3: Tikhonov-Regularized Operator Reconstruction}
		\STATE Evaluate basis functions at target $\boldsymbol{x}$: $\boldsymbol{b} = [u_1(\boldsymbol{x}), \dots, u_{M}(\boldsymbol{x})] \in \mathbb{C}^{1 \times M}$.
		\STATE Compute the optimal basis operator: $\boldsymbol{a}_* = \boldsymbol{b} \boldsymbol{V}^* (\boldsymbol{V} \boldsymbol{V}^* + \alpha \boldsymbol{I})^{-1}$, where $\boldsymbol{I} \in \mathbb{R}^{2N \times 2N}$ is the identity matrix.
		
		\STATE \COMMENT{Step 4: Solution Reconstruction}
		\STATE Assemble $\boldsymbol{h} = [f(\boldsymbol{x}_1), g(\boldsymbol{x}_1), \dots, f(\boldsymbol{x}_N), g(\boldsymbol{x}_N)]^T$.
		\RETURN $u^{(c)}(\boldsymbol{x}) = \mathrm{Re}(\boldsymbol{a}_* \boldsymbol{h})$.
	\end{algorithmic}
\end{algorithm}

\section{Theoretical Foundations and Convergence Analysis}
\label{sec:theory}
This section develops the analytical foundation of the proposed LbNM. For a biharmonic solution $u$, the Almansi decomposition $u= \phi + r^2\psi$ reduces the approximation problem to two harmonic components. The analysis therefore proceeds in three steps:

\begin{itemize}
    \item \textbf{Harmonic approximation}: We estimate the truncation error of normalized  harmonic polynomials on the fictitious disk $O_\rho$. 
    \item \textbf{Stable Almansi decomposition}: We control the two harmonic components in terms of the biharmonic solution.
    \item \textbf{Global convergence}: We combine the approximation estimates, boundary stability, collocation error, and Tikhonov regularization.
\end{itemize}

\subsection{Approximation Properties of Harmonic Bases}
Since the components $\phi$ and $\psi$ in the Almansi decomposition \eqref{eq:almansi_new} are  harmonic, we first establish the approximation estimates for harmonic polynomials on $O_\rho$.

We use the following stability estimate to transfer boundary
approximation error to the interior \cite{evans2022partial}:

\begin{lemma}[Interior Stability Estimate]
\label{lem: interior_stability_estimate}
Let $\Omega$ be a bounded Lipschitz domain. If $u \in H^1(\Omega)$ satisfies 
$\Delta u  = 0$ in $\Omega$, then
\begin{equation*}
\|u\|_{L^2(\Omega)} \leq C_{\text{stab}} \|u\|_{L^2(\partial \Omega)},
\end{equation*}
where $C_{\text{stab}}$ is a constant depending only on the geometry of $\Omega$.
\end{lemma}

Recall the form of the basis functions defined on the boundary of $O_\rho$ as follows:
\begin{equation*}
\varphi_n(r,\theta)
=
\frac{r^{|n|}e^{\mathrm{i}n\theta}}{\left\|r^{|n|}e^{\mathrm{i}n\theta}\right\|_{L^2(\partial O_\rho)}}
=
\frac{r^{|n|}e^{\mathrm{i}n\theta}}{\sqrt{2\pi\rho}\rho^{|n|}}.
\end{equation*}

The normalized basis has the following approximation property:

\begin{lemma}
\label{lem: harmonic_polynomial_approximation}
Let $0 < \rho < d$, let $O_\rho$ and $O_d$ be concentric disks centered at the origin. Assume $u\in H^1(O_d)$ be harmonic in $O_d$. Then, $u$ satisfies
\begin{equation*}
    u
    =
    \sum_{n\in \mathbb{Z}}\mu_n\varphi_n, \quad \text{with} \quad \|\boldsymbol{\mu}\|_{\ell^2} \leq \sqrt{\frac{\rho}{d}}\|u\|_{L^2(\partial O_d)}.
\end{equation*}
Furthermore, truncating the approximation to $L = 2K + 1$ terms, denoted by $u^{(L)}$, we have the error estimate
\begin{equation*}
    \|u - u^{(L)}\|_{L^2(O_\rho)}
    \leq
    C\|u\|_{L^2(\partial O_d)}\Big(\frac{\rho}{d}\Big)^{\frac{L}{2}},
\end{equation*}
where $C$ is a constant independent of $L$.
\end{lemma}

\begin{proof}
	For a harmonic function $u$, the trace of $u$ on $\partial O_d$ has the Fourier expansion
	\begin{equation*}
		u(d,\theta)
		=
		\sum_{n\in \mathbb{Z}}\hat{u}_n e^{\mathrm{i} n\theta}, \quad
		\hat{u}_n 
		=
		\frac{1}{2\pi}\int_0^{2\pi}u(d,\theta)e^{-\mathrm{i}n\theta}\mathrm{d}\theta,
		\end{equation*}
	By Parseval's identity, we have
	\begin{equation*}
		\|u\|_{L^2(\partial O_d)}^2
		=
		2\pi d\sum_{n\in \mathbb{Z}}|\hat{u}_n|^2.
		\end{equation*}
	The interior values of a harmonic function are given by the Poisson integral, which is equivalent to
	\begin{equation}
		u(r,\theta)
		=
		\sum_{n\in \mathbb{Z}}\hat{u}_n\Big(\frac{r}{d}\Big)^{|n|}e^{\mathrm{i} n\theta}, \quad
		0\leq r \leq d,
		\end{equation}
	Rewriting this as a linear combination of the normalized basis functions yields
	\begin{equation}
		u(r,\theta)
		=
			\sum_{n\in \mathbb{Z}}\mu_n\varphi_n(r,\theta),\quad
			0\leq r \leq \rho,
		\end{equation}
	where
	\begin{equation}
		\mu_n
		=
		\sqrt{2\pi\rho}\,\hat{u}_n\Big(\frac{\rho}{d}\Big)^{|n|},
		\end{equation}
	and
	\begin{equation}
		\begin{aligned}
		\|\boldsymbol{\mu}\|_{\ell^2}^2
		=
		\sum_{n\in\mathbb{Z}}|\mu_n|^2
		=
		2\pi\rho 	\sum_{n\in\mathbb{Z}}|\hat{u}_n|^2\Big(\frac{\rho}{d}\Big)^{2|n|}
		\leq
	2\pi\rho 	\sum_{n\in\mathbb{Z}}|\hat{u}_n|^2 = \frac{\rho}{d}	\|u\|_{L^2(\partial O_d)}^2.
		\end{aligned}
		\end{equation}	
		Choosing a truncation parameter $K \in \mathbb{N}_+$, we let  $L = 2K +1$ denote the number of modes. Thus
		\begin{equation*}
			u^{(L)}(r,\theta)
			=
			\sum_{n = -K}^{K}\mu_n \varphi_n(r,\theta),
			\end{equation*}
		where the coefficients $\boldsymbol{\mu}^{(L)} = (\mu_{-K}, \cdots, \mu_{K})$ are the truncation of $\boldsymbol{\mu}$. Obviously, we have
		\begin{equation*}
			 \|\boldsymbol{\mu}^{(L)}\|_{\ell^2} 
			 \leq
			 \|\boldsymbol{\mu}\|_{\ell^2} \leq \sqrt{\frac{\rho}{d}}\|u\|_{L^2(\partial O_d)}.
			\end{equation*}
		On $\partial O_\rho$,
		 $$(u - u^{(L)})(\rho,\theta) = \sum_{|n|> K}\mu_n \varphi_n(\rho,\theta),$$ 
		and hence we obtain the error estimate
		\begin{equation*}
			\begin{aligned}
				\|u - u^{(L)}\|_{L^2(\partial O_\rho)}^2
			&	=
				\sum_{|n|> K}|\mu_n|^2
				=
				2\pi\rho \sum_{|n|> K}|\hat{u}_n|^2\Big(\frac{\rho}{d}\Big)^{2|n|}\\
				&\leq
				2\pi\rho \Big(\frac{\rho}{d}\Big)^{2K} \sum_{|n|> K}|\hat{u}_n|^2
				\leq
			\|u\|_{L^2(\partial O_d)}^2\Big(\frac{\rho}{d}\Big)^{2K+1}.
				\end{aligned}
			\end{equation*}
The claimed interior estimate now follows Lemma \ref{lem: interior_stability_estimate}.

	\end{proof}

We also recall the quantitative Runge approximation theorem for harmonic functions from \cite{ruland2019quantitative}. 

\begin{lemma}[
Quantitative Runge approximation,
specialized from~\cite{ruland2019quantitative}
]
\label{lem:quantitative-runge-approximation}

Let $D_1\Subset D_2$ be bounded Lipschitz domains such that
$D_2\setminus\overline{D_1}$ is connected, and let
$\Gamma\subset\partial D_2$ be nonempty and relatively open. Define
\begin{equation*}
S_1
\coloneqq
\left\{
w\in H^1(D_1):
\Delta w=0 \text{ in }D_1
\right\},
\end{equation*}
and
\begin{equation*}
S_2
\coloneqq
\left\{
w\in H^1(D_2):
\Delta w=0 \text{ in }D_2,\ 
\operatorname{supp}(w|_{\partial D_2})
\subset\overline{\Gamma}
\right\}.
\end{equation*}

There exist constants $C>1$ and $\theta>0$, depending on
$D_1$, $D_2$, and $\Gamma$, such that for every
$\varepsilon\in(0,1)$ and every $h\in S_1$, there exists
$w\in S_2$ satisfying
\begin{equation}
\label{eq:runge-exponential}
\|h-w\|_{L^2(D_1)}
\leq
\varepsilon\|h\|_{H^1(D_1)},
\qquad
\|w\|_{H^{1/2}(\partial D_2)}
\leq
C e^{C\varepsilon^{-\theta}}
\|h\|_{L^2(D_1)}.
\end{equation}

If, in addition,
$D_1\Subset\widehat D\Subset D_2$, with $\widehat D$ Lipschitz,
and $\widehat h\in H^1(\widehat D)$ is harmonic in $\widehat D$,
then there exist constants $C,\beta\geq1$, depending also on
$\widehat D$, and $w\in S_2$ such that
\begin{equation}
\label{eq:runge-polynomial}
\|\widehat h|_{D_1}-w\|_{L^2(D_1)}
\leq
\varepsilon\|\widehat h\|_{H^1(\widehat D)},
\qquad
\|w\|_{H^{1/2}(\partial D_2)}
\leq
C\varepsilon^{-\beta}
\|\widehat h\|_{L^2(D_1)}.
\end{equation}
\end{lemma}
    
% \begin{lemma}[Quantitative Runge Approximation] 
% \label{lem: quantitative_runge_approximation}
% Let $\Omega \subset \subset \bar{\Omega} \subset \subset \Omega_{ext}$ be bounded open Lipschitz sets. Define the spaces of solutions:
% \begin{equation*}
% S_1 \coloneqq \{w \in H^1(\Omega) : \Delta w = 0 \text{ in } \Omega\},
% \end{equation*}
% \begin{equation*}
% S_2 \coloneqq \{h \in H^1(\Omega_{ext}) : \Delta h = 0 \text{ in } \Omega_{ext}\}.
% \end{equation*}
% Then, for any $\epsilon \in (0,1)$, there exists a constant $C$, parameters $\theta > 0$ and $\beta \geq 1$ such that:

% \begin{enumerate}
% \item For any $h \in S_1$, there exists $u \in S_2$ satisfying:
% \begin{equation*}
% \|h - u\|_{L^2(\Omega)} \leq \epsilon \|h\|_{H^1(\Omega)},\quad
% \|u\|_{H^{1/2}(\partial \Omega_{ext})} \leq C e^{C \epsilon^{-\theta}} \|h\|_{L^2(\Omega)}.
% \end{equation*}

% \item If $h$ can be continued to $\tilde{h} \in H^1(\tilde{\Omega})$ satisfying the harmonic equation, then there exists $u \in S_2$ satisfying a polynomial bound:
% \begin{equation*}
% \|\tilde{h} - u\|_{L^2(\Omega)} \leq \epsilon \|\tilde{h}\|_{H^1(\Omega)},\quad
% \|u\|_{H^{1/2}(\partial \Omega_{ext})} \leq C \epsilon^{-\beta} \|\tilde{h}\|_{L^2(\Omega)}.
% \end{equation*}
% \end{enumerate}

% \end{lemma}

\subsection{Almansi Decomposition}
We next transfer the biharmonic estimates to harmonic functions through the Almansi decomposition \cite{almansi1899sull}. For completeness, we state and prove the Sobolev stability result needed below.

Following \cite{CostabelDauge2015}, we use the following definition.

\begin{definition}
\label{def:strict-star-shaped}
Let $\Omega\subset\mathbb{R}^d$ be a domain. We say that $\Omega$ is
\emph{strictly star-shaped} if there exists an open ball
$B\subset\Omega$ such that
\begin{equation*}
    [\boldsymbol{y},\boldsymbol{x}]\subset\Omega
    \qquad
    \text{for every }\boldsymbol{y}\in B
    \text{ and }\boldsymbol{x}\in\Omega,
\end{equation*}
where
\begin{equation*}
    [\boldsymbol{y},\boldsymbol{x}]
    :=
    \bigl\{(1-t)\boldsymbol{y}+t\boldsymbol{x}:\ 0\leq t\leq1\bigr\}
\end{equation*}
denotes the closed line segment joining $\boldsymbol{y}$ and $\boldsymbol{x}$.

We say that $\Omega$ is \emph{strictly star-shaped with respect to the origin} if the ball $B$ can be chosen as $ B=B_a(0)$ for some $a>0$.
\end{definition}

\begin{lemma}
\label{lem:strict-star-shaped-dilation}
Let $\Omega\subset\mathbb{R}^d$ be a bounded domain that is strictly star-shaped with respect to the origin. 
Then, for every $\rho\in(0,1)$,
\begin{equation*}
    \operatorname{dist}
    \bigl(\overline{\rho\Omega},\Omega^c\bigr)
    \geq
    (1-\rho)a.
\end{equation*}
In particular,
$\overline{\rho\Omega}\Subset\Omega.$

\end{lemma}

\begin{proof}
Fix $\rho\in(0,1)$ and $\boldsymbol{x}\in\Omega$. We claim that 
\begin{equation*}
    B_{(1-\rho)a}(\rho \boldsymbol{x})\subset\Omega.
\end{equation*}
Indeed, let $\boldsymbol{z}\in B_{(1-\rho)a}(\rho \boldsymbol{x})$. 
Then $|\boldsymbol{z}-\rho \boldsymbol{x}|<(1-\rho)a$.
Define $ \boldsymbol{y}=\frac{\boldsymbol{z}-\rho \boldsymbol{x}}{1-\rho}.$
It follows that
$|\boldsymbol{y}|<a,$ and hence $\boldsymbol{y}\in B_a(0)$. Moreover,
$\boldsymbol{z} =\rho \boldsymbol{x}+(1-\rho)\boldsymbol{y}$.
Thus $\boldsymbol{z}$ lies on the line segment joining $\boldsymbol{y}$ and $\boldsymbol{x}$. Since $\Omega$ is star-shaped with respect to every point of $B_a(0)$, we obtain
$\boldsymbol{z}\in[\boldsymbol{y},\boldsymbol{x}]\subset\Omega.$
This proves
\begin{equation*}
    B_{(1-\rho)a}(\rho \boldsymbol{x})\subset\Omega
    \qquad
    \text{for every }\boldsymbol{x}\in\Omega.
\end{equation*}

Consequently,
\begin{equation*}
    \operatorname{dist}(\rho \boldsymbol{x},\Omega^c)
    \geq
    (1-\rho)a
    \qquad
    \text{for every }\boldsymbol{x}\in\Omega,
\end{equation*}
and therefore
\begin{equation*}
    \operatorname{dist}(\rho\Omega,\Omega^c)
    \geq
    (1-\rho)a.
\end{equation*}
Since the function $ \boldsymbol{z}\mapsto\operatorname{dist}(\boldsymbol{z},\Omega^c)$
is continuous, the same estimate holds on the closure:
\begin{equation*}
    \operatorname{dist}
    \bigl(\overline{\rho\Omega},\Omega^c\bigr)
    \geq
    (1-\rho)a>0.
\end{equation*}
In particular, $\overline{\rho\Omega}\subset\Omega.$
Finally, since $\Omega$ is bounded, $\overline{\rho\Omega}$ is compact. Hence $ \overline{\rho\Omega}\Subset\Omega.$

\end{proof}

With this geometric property, we now establish a crucial gradient estimate for harmonic functions
\begin{lemma}
\label{lem:euler-harmonic-estimate}
Let $\Omega\subset\mathbb{R}^2$ be a bounded Lipschitz domain that is
strictly star-shaped with respect to the origin.  Let
\[
    v\in L^2(\Omega;\mathbb{R})
\]
satisfy
\[
    \Delta v=0
    \qquad\text{in }\mathcal{D}'(\Omega)
\]
and
\[
    \boldsymbol{x}\cdot\nabla v\in L^2(\Omega).
\]
Then $v\in H^1(\Omega)$ and
\begin{equation}
\label{eq:euler-harmonic-estimate}
    \|v\|_{H^1(\Omega)}
    \leq
    C_{\Omega}
    \left(
        \|v\|_{L^2(\Omega)}
        +
        \|\boldsymbol{x}\cdot\nabla v\|_{L^2(\Omega)}
    \right),
\end{equation}
where $C_{\Omega}>0$ depends only on $\Omega$.
\end{lemma}

\begin{proof}
Since $v\in L^2(\Omega)$, by Weyl's lemma \cite[Corollary 2.2.1]{JostPDE}, we henceforth regard $v$ as a smooth harmonic function in $\Omega$. For $\boldsymbol{x}=(x_1,x_2)$, denote
$   z=x_1+\mathrm{i}x_2$,
and define
\begin{equation*}
    g(z)=v_{x_1}-\mathrm{i}v_{x_2}.
\end{equation*}
A direct calculation gives
\begin{equation*}
\begin{aligned}
    \partial_{\bar z}g
    =
    \frac12
    \bigl(\partial_{x_1}+\mathrm{i}\partial_{x_2}\bigr)
    \bigl(v_{x_1}-\mathrm{i}v_{x_2}\bigr) 
    =
    \frac12
    \bigl(v_{x_1x_1}+v_{x_2x_2}\bigr)
    =
    \frac12\Delta v
    =
    0.
\end{aligned}
\end{equation*}
Hence $g$ is analytic in $\Omega$, and therefore
$   F(z)=zg(z) $
is analytic in $\Omega$. Writing
\begin{equation*}
    F=q+\mathrm{i}p,
\end{equation*}
we obtain
\begin{equation*}
\begin{aligned}
    F(z)
    =
    (x_1+\mathrm{i}x_2)
    (v_{x_1}-\mathrm{i}v_{x_2}) 
    =
    (x_1v_{x_1}+x_2v_{x_2})
    +
    \mathrm{i}(x_2v_{x_1}-x_1v_{x_2}).
\end{aligned}
\end{equation*}
Consequently,
\begin{equation}
\label{eq:pq-definitions}
    q=\boldsymbol{x}\cdot\nabla v,
    \qquad
    p=x_2v_{x_1}-x_1v_{x_2},
\end{equation}
 Weyl's lemma implies that $p$ and $q$ are smooth in
$\Omega$, and hence $p,q\in L^2_{\mathrm{loc}}(\Omega)$. To apply the Friedrichs--Velte inequality \cite{CostabelDauge2015,Costabel2017}  and prove that $p\in L^2(\Omega)$,  we work on compactly contained dilations of $\Omega$. 
Fix $\rho_0\in(0,1)$ and, for $\rho\in[\rho_0,1)$, set $\Omega_\rho=\rho\Omega. $ 
Since $\overline{\Omega_\rho}\Subset\Omega$ by Lemma \ref{lem:strict-star-shaped-dilation}, one has
$p,q\in L^2(\Omega_\rho)$.
Let
\begin{equation*}
    p_\rho
    =
    \frac{1}{|\Omega_\rho|}
    \int_{\Omega_\rho}p\, \mathrm{d} \boldsymbol{x},
    \quad
    q_\rho
    =
    \frac{1}{|\Omega_\rho|}
    \int_{\Omega_\rho}q\,\mathrm{d}\boldsymbol{x}.
\end{equation*}
The function
\[
\begin{aligned}
    -\mathrm{i}
    \bigl(F-(q_\rho+\mathrm{i}p_\rho)\bigr)
    &=
    (p-p_\rho)-\mathrm{i}(q-q_\rho)
\end{aligned}
\]
is analytic in $\Omega_\rho$, and
\[
    \int_{\Omega_\rho}(p-p_\rho)\,\mathrm{d}\boldsymbol{x}=0.
\]
The Friedrichs--Velte inequality for conjugate harmonic functions
therefore yields
\begin{equation}
\label{eq:friedrichs-velte-rho}
    \|p-p_\rho\|_{L^2(\Omega_\rho)}
    \leq
    C_F(\Omega_\rho)
    \|q-q_\rho\|_{L^2(\Omega_\rho)}.
\end{equation}
 The Friedrichs--Velte constant $C_F(\Omega_\rho)$ is invariant under dilations. To see
this explicitly, let \(H=h+\mathrm{i}k\) be analytic in
\(\Omega_\rho=\rho\Omega\), and define the pullbacks
\begin{equation*}
    \widetilde h(\boldsymbol{y})=h(\rho \boldsymbol{y}),
    \qquad
    \widetilde k(\boldsymbol{y})=k(\rho \boldsymbol{y}),
    \quad \boldsymbol{y}\in\Omega.
\end{equation*}
Then
\begin{equation*}
    \widetilde H(\boldsymbol{y})
    =
    H(\rho \boldsymbol{y})
    =
    \widetilde h(\boldsymbol{y})+\mathrm{i}\widetilde k(\boldsymbol{y})
\end{equation*}
is analytic in \(\Omega\).

The corresponding mean values are preserved by this pullback.
Indeed, using \(\boldsymbol{x}=\rho \boldsymbol{y}\) and \(\mathrm{d}\boldsymbol{x}=\rho^2\,\mathrm{d}\boldsymbol{y}\), we obtain
\begin{equation*}
\begin{aligned}
    h_{\Omega_\rho}
    :=
    \frac{1}{|\Omega_\rho|}
    \int_{\Omega_\rho}h(\boldsymbol{x})\,\mathrm{d}\boldsymbol{x} 
    =
    \frac{1}{\rho^2|\Omega|}
    \int_\Omega h(\rho \boldsymbol{y} )\rho^2\,\mathrm{d}\boldsymbol{y} 
    =
    \frac{1}{|\Omega|}
    \int_\Omega\widetilde h(\boldsymbol{y} )\,\mathrm{d}\boldsymbol{y}
    =
    \widetilde h_\Omega.
\end{aligned}
\end{equation*}
Similarly,
$ k_{\Omega_\rho}=\widetilde k_\Omega$. 
Consequently, the  \(L^2\) norms satisfy
\begin{equation*}
\begin{aligned}
    \|h-h_{\Omega_\rho}\|_{L^2(\Omega_\rho)}^2
    =
    \int_{\Omega_\rho}
    |h(\boldsymbol{x} )-h_{\Omega_\rho}|^2\,\mathrm{d}\boldsymbol{x} 
    =
    \rho^2
    \int_\Omega
    |\widetilde h(\boldsymbol{y})-\widetilde h_\Omega|^2\,\mathrm{d}\boldsymbol{y}
    =
    \rho^2
    \|\widetilde h-\widetilde h_\Omega\|_{L^2(\Omega)}^2,
\end{aligned}
\end{equation*}
The same calculation gives
\[
    \|k-k_{\Omega_\rho}\|_{L^2(\Omega_\rho)}^2
    =
    \rho^2
    \|\widetilde k-\widetilde k_\Omega\|_{L^2(\Omega)}^2.
\]

Applying the Friedrichs--Velte inequality on \(\Omega\) to
\(\widetilde H\), we find
\begin{equation*}
\begin{aligned}
    \|h-h_{\Omega_\rho}\|_{L^2(\Omega_\rho)}
    =
    \rho
    \|\widetilde h-\widetilde h_\Omega\|_{L^2(\Omega)} 
    \leq
    \rho C_F(\Omega)
    \|\widetilde k-\widetilde k_\Omega\|_{L^2(\Omega)} 
    =
    C_F(\Omega)
    \|k-k_{\Omega_\rho}\|_{L^2(\Omega_\rho)}.
\end{aligned}
\end{equation*}
Since this holds for every admissible analytic function on
\(\Omega_\rho\), the optimality of the Friedrichs--Velte constant
implies
\[
    C_F(\Omega_\rho)\leq C_F(\Omega).
\]

Conversely, starting with an analytic function on \(\Omega\) and
rescaling it to \(\Omega_\rho\), the same argument with dilation factor
\(1/\rho\) gives
\[
    C_F(\Omega)\leq C_F(\Omega_\rho).
\]
Therefore,
\[
    C_F(\Omega_\rho)=C_F(\Omega),
\]
and the constant in \eqref{eq:friedrichs-velte-rho} is independent of
\(\rho\).

We next control the mean value $p_\rho$. Since the origin is an
interior point of $\Omega$, there exists $R>0$ such that
$\overline{B_R(0)}\subset\Omega$.
For every $\rho\in[\rho_0,1)$,
\[
    B_{\rho_0R}(0)
    \subset
    B_{\rho R}(0)
    =
    \rho B_R(0)
    \subset
    \Omega_\rho.
\]
Set
$r_*=\rho_0R$.
Because $p-p_\rho$ is harmonic in $\Omega_\rho$, the mean-value
property and the Cauchy--Schwarz inequality give
\begin{equation*}
\begin{aligned}
    |(p-p_\rho)(0)|
    =
    \Big|
        \frac{1}{|B_{r_*}|}
        \int_{B_{r_*}}(p-p_\rho)\,\mathrm{d}\boldsymbol{x}
    \Big| 
    \leq
    |B_{r_*}|^{-1/2}
    \|p-p_\rho\|_{L^2(B_{r_*})} 
    \leq
    C_\Omega
    \|p-p_\rho\|_{L^2(\Omega_\rho)}.
\end{aligned}
\end{equation*}
On the other hand, $F(0)=0$, and hence $p(0)=0$. Therefore,
$(p-p_\rho)(0)=-p_\rho,$ so that
\begin{equation}
\label{eq:p-mean-bound}
    |p_\rho|
    \leq
    C_\Omega
    \|p-p_\rho\|_{L^2(\Omega_\rho)}.
\end{equation}
Using \eqref{eq:p-mean-bound}, we obtain
\begin{equation*}
\begin{aligned}
    \|p\|_{L^2(\Omega_\rho)}
    \leq
    \|p-p_\rho\|_{L^2(\Omega_\rho)}
    +
    |\Omega_\rho|^{1/2}|p_\rho| 
    \leq
    C_\Omega
    \|p-p_\rho\|_{L^2(\Omega_\rho)}
    \leq
    C_\Omega
    \|q-q_\rho\|_{L^2(\Omega_\rho)}.
\end{aligned}
\end{equation*}
Since subtraction of the mean is the $L^2$-orthogonal projection onto the zero-mean subspace,
\[
    \|q-q_\rho\|_{L^2(\Omega_\rho)}
    \leq
    \|q\|_{L^2(\Omega_\rho)}.
\]
It follows that
\begin{equation}
\label{eq:p-by-q-rho}
    \|p\|_{L^2(\Omega_\rho)}
    \leq
    C_\Omega
    \|q\|_{L^2(\Omega_\rho)},
\end{equation}
with a constant independent of
$\rho\in[\rho_0,1)$. 
The family
$\{\Omega_\rho\}_{\rho\in[\rho_0,1)}$
is increasing and exhausts $\Omega$ as $\rho\to1^-$. Hence, by the
monotone convergence theorem and \eqref{eq:p-by-q-rho},
\begin{equation*}
\begin{aligned}
    \int_\Omega |p|^2\,\mathrm{d}\boldsymbol{x}
    =
    \lim_{\rho\to1^-}
    \int_{\Omega_\rho}|p|^2\,\mathrm{d}\boldsymbol{x}
    \leq
    C_\Omega
    \lim_{\rho\to1^-}
    \int_{\Omega_\rho}|q|^2\,\mathrm{d}\boldsymbol{x} 
    =
    C_\Omega
    \int_\Omega|q|^2\,\mathrm{d}\boldsymbol{x}.
\end{aligned}
\end{equation*}
In particular, $p\in L^2(\Omega)$ and
\begin{equation}
\label{eq:tangential-by-radial}
    \|p\|_{L^2(\Omega)}
    \leq
    C_\Omega\|q\|_{L^2(\Omega)}.
\end{equation}

The pointwise identity
\begin{equation*}
\begin{aligned}
    (x_1v_{x_1}+x_2v_{x_2})^2
    +(x_2v_{x_1}-x_1v_{x_2})^2 
    =
    |\boldsymbol{x}|^2|\nabla v|^2
\end{aligned}
\end{equation*}
and \eqref{eq:pq-definitions}--\eqref{eq:tangential-by-radial} imply
\begin{equation*}
\begin{aligned}
    \||\boldsymbol{x}|\nabla v\|_{L^2(\Omega)}^2
    =
    \|q\|_{L^2(\Omega)}^2
    +
    \|p\|_{L^2(\Omega)}^2 
    \leq
    C_\Omega
    \|q\|_{L^2(\Omega)}^2.
\end{aligned}
\end{equation*}
Thus
\begin{equation}
\label{eq:weighted-gradient-estimate}
    \||\boldsymbol{x}|\nabla v\|_{L^2(\Omega)}
    \leq
    C_\Omega
    \|\boldsymbol{x}\cdot\nabla v\|_{L^2(\Omega)}.
\end{equation}

Finally, choose $r>0$ such that $\overline{B_{2r}(0)}\subset\Omega.$
The Caccioppoli inequality \cite[Lemma 14.3.2]{JostPDE} for harmonic functions gives
\begin{equation*}
\begin{aligned}
    \|\nabla v\|_{L^2(B_r)}
    \leq
    \frac{C}{r}
    \|v\|_{L^2(B_{2r})} 
    \leq
    C_\Omega\|v\|_{L^2(\Omega)};
\end{aligned}
\end{equation*}
see, for example, \cite[Chapter~8]{GilbargTrudinger2001}.
For $\boldsymbol{x}\in\Omega\setminus B_r(0)$, one has
$|\boldsymbol{x}|\geq r$, and hence
\begin{equation*}
\begin{aligned}
    \|\nabla v\|_{L^2(\Omega\setminus B_r)}
    \leq
    \frac{1}{r}
    \||\boldsymbol{x}|\nabla v\|_{L^2(\Omega)} 
    \leq
    C_\Omega
    \|\boldsymbol{x}\cdot\nabla v\|_{L^2(\Omega)}.
\end{aligned}
\end{equation*}
Combining the estimates on $B_r$ and
$\Omega\setminus B_r$ yields
\begin{equation*}
    \|\nabla v\|_{L^2(\Omega)}
    \leq
    C_\Omega
    \left(
        \|v\|_{L^2(\Omega)}
        +
        \|\boldsymbol{x}\cdot\nabla v\|_{L^2(\Omega)}
    \right).
\end{equation*}
Therefore the estimate \eqref{eq:euler-harmonic-estimate} follows.
\end{proof}

Using this estimate, we now prove the Sobolev stability of the Almansi decomposition.

\begin{theorem}[Sobolev-stable Almansi decomposition]
\label{thm:stable-almansi-decomposition}
Let $\Omega\subset\mathbb{R}^2$ satisfy the assumptions of
Lemma~\ref{lem:euler-harmonic-estimate}. Suppose that
$u\in H^2(\Omega)$ and
\begin{equation*}
    \Delta^2u=0
    \qquad\text{in }\mathcal{D}'(\Omega).
\end{equation*}
Then there exists a unique pair of harmonic functions $\phi,\psi\in H^1(\Omega)$ such that
\begin{equation}
\label{eq:almansi-decomposition}
    u=\phi+|\boldsymbol{x}|^2\psi
    \qquad\text{a.e. in }\Omega.
\end{equation}
Moreover,
\begin{equation}
\label{eq:stable-almansi-estimate}
    \|\phi\|_{H^1(\Omega)}
    +
    \|\psi\|_{H^1(\Omega)}
    \leq
    C_\Omega\|u\|_{H^2(\Omega)}.
\end{equation}
\end{theorem}

\begin{proof}
Set $ w = \Delta u$.
Since $u\in H^2(\Omega)$, one has $w\in L^2(\Omega)$. Moreover,
\[
    \Delta w=\Delta^2u=0
    \qquad\text{in }\mathcal{D}'(\Omega).
\]
By Weyl's lemma, $ w\in C^\infty(\Omega)$. 
Define
\begin{equation}
\label{eq:psi-definition}
    \psi(\boldsymbol{x})
    =
    \frac14\int_0^1w(t\boldsymbol{x})\,\mathrm{d}t.
\end{equation}
The definition is meaningful because $\Omega$ is star-shaped with
respect to the origin, so that
\begin{equation*}
    t\boldsymbol{x}\in\Omega
    \qquad
    \text{for all }\boldsymbol{x}\in\Omega
    \text{ and }t\in[0,1].
\end{equation*}

We first verify that $\psi$ is harmonic. Since $w$ is smooth in
$\Omega$, differentiation under the integral sign is valid locally,
and
\begin{equation*}
\begin{aligned}
    \Delta\psi(\boldsymbol{x})
    =
    \frac14
    \int_0^1
    t^2(\Delta w)(t\boldsymbol{x})\,\mathrm{d}t 
    =0.
\end{aligned}
\end{equation*}
Furthermore,
\begin{equation*}
\begin{aligned}
    \psi(\boldsymbol{x})
    +
    \boldsymbol{x}\cdot\nabla\psi(\boldsymbol{x})
    =
    \frac14
    \int_0^1
    \left[
        w(t\boldsymbol{x})
        +
        t\boldsymbol{x}\cdot
        \nabla w(t\boldsymbol{x})
    \right]\mathrm{d}t 
    =
    \frac14
    \int_0^1
    \frac{\mathrm{d}}{\mathrm{d}t}
    \left[
        t\,w(t\boldsymbol{x})
    \right]\mathrm{d}t 
    =
    \frac14w(\boldsymbol{x}).
\end{aligned}
\end{equation*}
Consequently, we have
\begin{equation}
\label{eq:euler-equation-psi}
    (I+\boldsymbol{x}\cdot\nabla)\psi
    =
    \frac14w
    =
    \frac14\Delta u.
\end{equation}

We next prove that $\psi\in L^2(\Omega)$.
Choose $R>0$ such that $\overline{B_{2R}(0)}\subset\Omega.$
Since $\Omega$ is bounded, there exists $a\in(0,1)$ such that
$a\Omega\subset B_R(0).$
For $0<t\leq a$, the interior estimate for harmonic functions gives
\begin{equation*}
\begin{aligned}
    \|w(t\,\cdot)\|_{L^2(\Omega)}
    \leq
    |\Omega|^{1/2}
    \sup_{B_R(0)}|w| 
    \leq
    C_\Omega
    \|w\|_{L^2(B_{2R}(0))} 
    \leq
    C_\Omega
    \|w\|_{L^2(\Omega)}.
\end{aligned}
\end{equation*}
For $a\leq t\leq1$, the change of variables
$\boldsymbol{y}=t\boldsymbol{x}$ gives
\begin{equation*}
\begin{aligned}
    \|w(t\,\cdot)\|_{L^2(\Omega)}^2
    =
    t^{-2}
    \int_{t\Omega}|w(\boldsymbol{y})|^2\,d\boldsymbol{y} 
    \leq
    a^{-2}
    \|w\|_{L^2(\Omega)}^2.
\end{aligned}
\end{equation*}
Splitting the integral at $t=a$, we obtain
\begin{equation}
\label{eq:psi-l2-bound}
\begin{aligned}
\|\psi\|_{L^2(\Omega)}
&\leq
\frac14\int_0^a
\|w(t\,\cdot)\|_{L^2(\Omega)}\,\mathrm{d}t
+
\frac14\int_a^1
\|w(t\,\cdot)\|_{L^2(\Omega)}\,\mathrm{d}t \\
&\leq
C_\Omega
\left(
a+\frac{1-a}{a}
\right)
\|w\|_{L^2(\Omega)}
\leq
C_\Omega\|w\|_{L^2(\Omega)}.
\end{aligned}
\end{equation}
Since $\psi$ is harmonic and
$\boldsymbol{x}\cdot\nabla\psi
    =
    \frac14w-\psi
    \in L^2(\Omega)$. 
Lemma \ref{lem:euler-harmonic-estimate},
\eqref{eq:euler-equation-psi}, and
\eqref{eq:psi-l2-bound} give
\begin{equation*}
\begin{aligned}
    \|\psi\|_{H^1(\Omega)}
    \leq
    C_\Omega
    \left(
        \|\psi\|_{L^2(\Omega)}
        +
        \|\boldsymbol{x}\cdot\nabla\psi\|_{L^2(\Omega)}
    \right) 
    \leq
    C_\Omega
    \left(
        \|\psi\|_{L^2(\Omega)}
        +
        \|w\|_{L^2(\Omega)}
    \right) 
    \leq
    C_\Omega\|w\|_{L^2(\Omega)}.
\end{aligned}
\end{equation*}
On the other hand, define
\begin{equation}
\label{eq:phi-definition}
    \phi
    =
    u-|\boldsymbol{x}|^2\psi.
\end{equation}
Since $\Omega$ is bounded and $\psi\in H^1(\Omega)$, multiplication
by $|\boldsymbol{x}|^2$ maps $H^1(\Omega)$ continuously into itself. Thus $\phi\in H^1(\Omega).$
For a harmonic function $\psi$ in two dimensions,
\begin{equation*}
\begin{aligned}
    \Delta\bigl(|\boldsymbol{x}|^2\psi\bigr)
    =
    (\Delta|\boldsymbol{x}|^2)\psi
    +
    2\nabla|\boldsymbol{x}|^2\cdot\nabla\psi
    +
    |\boldsymbol{x}|^2\Delta\psi 
    =
    4\psi
    +
    4\boldsymbol{x}\cdot\nabla\psi 
    =
    4(I+\boldsymbol{x}\cdot\nabla)\psi.
\end{aligned}
\end{equation*}
Using \eqref{eq:euler-equation-psi}, we obtain
\begin{equation*}
    \Delta\bigl(|\boldsymbol{x}|^2\psi\bigr)
    =
    w
    =
    \Delta u
\end{equation*}
in $\mathcal{D}'(\Omega)$. It follows from
\eqref{eq:phi-definition} that
\begin{equation*}
    \Delta\phi=0
    \qquad\text{in }\mathcal{D}'(\Omega).
\end{equation*}
Hence $\phi$ is harmonic, and
\eqref{eq:almansi-decomposition} holds.
The estimate for $\psi$ gives
\begin{equation*}
    \|\psi\|_{H^1(\Omega)}
    \leq
    C_\Omega
    \|\Delta u\|_{L^2(\Omega)}
    \leq
    C_\Omega
    \|u\|_{H^2(\Omega)}.
\end{equation*}
Moreover,
\begin{equation*}
\begin{aligned}
    \|\phi\|_{H^1(\Omega)}
    \leq
    \|u\|_{H^1(\Omega)}
    +
    \||\boldsymbol{x}|^2\psi\|_{H^1(\Omega)} 
    \leq
    C_\Omega
    \left(
        \|u\|_{H^1(\Omega)}
        +
        \|\psi\|_{H^1(\Omega)}
    \right) 
    &\leq
    C_\Omega\|u\|_{H^2(\Omega)}.
\end{aligned}
\end{equation*}
This proves \eqref{eq:stable-almansi-estimate}.

It remains to prove uniqueness. Suppose that
\begin{equation*}
    u
    =
    \phi_1+|\boldsymbol{x}|^2\psi_1
    =
    \phi_2+|\boldsymbol{x}|^2\psi_2,
\end{equation*}
where $\phi_j,\psi_j\in H^1(\Omega)$ are harmonic. Set
\begin{equation*}
    \Phi=\phi_1-\phi_2,
    \qquad
    \Psi=\psi_1-\psi_2.
\end{equation*}
Then
\begin{equation*}
    \Phi+|\boldsymbol{x}|^2\Psi=0.
\end{equation*}
Taking the Laplacian in the sense of distributions gives
\begin{equation*}
    4(I+\boldsymbol{x}\cdot\nabla)\Psi=0.
\end{equation*}
Since $\Psi$ is harmonic, it is smooth in $\Omega$. For fixed
$\boldsymbol{x}\in\Omega$ and $t\in(0,1]$, we have
\begin{equation*}
\begin{aligned}
    \frac{\mathrm{d}}{\mathrm{d}t}
    \left[
        t\Psi(t\boldsymbol{x})
    \right]
    =
    \Psi(t\boldsymbol{x})
    +
    t\boldsymbol{x}\cdot
    \nabla\Psi(t\boldsymbol{x}) 
    =0.
\end{aligned}
\end{equation*}
Thus $t\Psi(t\boldsymbol{x})$ is independent of $t$. Since $\Psi$ is bounded in a neighborhood of the origin,
\begin{equation*}
    \lim_{t\to0^+}t\Psi(t\boldsymbol{x})=0.
\end{equation*}
Therefore
\begin{equation*}
    \Psi(\boldsymbol{x})=0
    \qquad
    \text{for every }\boldsymbol{x}\in\Omega.
\end{equation*}
It follows that $\Psi=0$ and then $\Phi=0$. Hence the decomposition is unique.
\end{proof}

\begin{remark}
\label{rem:dimension-two}
The proof of Lemma~\ref{lem:euler-harmonic-estimate} uses the
two-dimensional complex-analytic representation of conjugate harmonic functions. Consequently, Theorem~\ref{thm:stable-almansi-decomposition}, in the form stated above, is a two-dimensional result. The classical Almansi decomposition exists in arbitrary dimension, but the present $H^1$-stability argument does not extend directly to dimensions $d\geq3$.
\end{remark}

\subsection{Convergence and Stability of Biharmonic LbNM}
The stability estimate in Theorem \ref{thm:stable-almansi-decomposition} controls the harmonic components in terms of the biharmonic solution. We now combine this estimate with harmonic approximation and regularized boundary fitting.

\begin{lemma}
\label{lem:biharmonic_polynomial_approximation}
Let $0 < \rho < d$, let $O_\rho$ and $O_d$ be concentric disks centered at the origin, and let $u \in H^2(O_d)$ be biharmonic in $O_d$, with Almansi decomposition $u = \phi + r^2 \psi$, where $\phi$ and $\psi$ are harmonic in $O_d$. Define the biharmonic basis functions by
\begin{equation*}
    \Phi_n^{(1)} = \varphi_n, \quad \Phi_n^{(2)} = r^2\varphi_n, \quad n \in \mathbb{Z}.
\end{equation*}
Then, $u$ can be represented as the following convergent series:
\begin{equation*}
    u = \sum_{n \in \mathbb{Z}} \left( \xi_n \Phi_n^{(1)} + \eta_n \Phi_n^{(2)} \right),
\end{equation*}
where the expansion coefficients $\boldsymbol{\xi} = \{\xi_n\}_{n \in \mathbb{Z}}$ and $\boldsymbol{\eta} = \{\eta_n\}_{n \in \mathbb{Z}}$ satisfy the bound
\begin{equation}
\label{eq:coefficient_bound}     
    \left\|(\boldsymbol{\xi}, \boldsymbol{\eta})\right\|_{\ell^2\times\ell^2} 
    \leq 
    \sqrt{\frac{\rho}{d}} \left( \|\phi\|_{L^2(\partial O_d)} + \|\psi\|_{L^2(\partial O_d)} \right).
\end{equation}
Furthermore, let $L= 2K+1$ be the number of Fourier modes retained in each harmonic component, so that the truncated biharmonic expansion $u^{(L)}$ contains $M = 2L$ basis functions.   Then
\begin{equation}
\label{eq:polynomial_approximation}
    \left\|u - u^{(L)}\right\|_{L^2(O_\rho)}
    \leq
    C \left( \|\phi\|_{L^2(\partial O_d)} + \|\psi\|_{L^2(\partial O_d)} \right) \left(\frac{\rho}{d}\right)^{\frac{L}{2}},
\end{equation}
where the constant $C > 0$ is independent of $K$ (and hence of $L$ and $M$).
\end{lemma}

\begin{proof}
By the assumed Almansi decomposition,
\begin{equation*}
    u(\boldsymbol{x}) = \phi(\boldsymbol{x}) + |\boldsymbol{x}|^2\psi(\boldsymbol{x}), 
\end{equation*}
where $\phi$ and $\psi$ are harmonic functions satisfying $\Delta\phi = 0$ and $\Delta\psi = 0$ in $O_d$. Applying Lemma \ref{lem: harmonic_polynomial_approximation} to the individual harmonic components $\phi$ and $\psi$, we obtain the following expansions:
\begin{equation*}
    \phi = \sum_{n \in \mathbb{Z}} \xi_n\varphi_n, \quad \psi = \sum_{n \in \mathbb{Z}} \eta_n\varphi_n,
\end{equation*}
with the coefficients satisfying the $\ell^2$-bounds:
\begin{equation*}
    \|\boldsymbol{\xi}\|_{\ell^2} 
    \leq \sqrt{\frac{\rho}{d}}\|\phi\|_{L^2(\partial O_d)}, \quad \|\boldsymbol{\eta}\|_{\ell^2} \leq \sqrt{\frac{\rho}{d}}\|\psi\|_{L^2(\partial O_d)},
\end{equation*}
and the corresponding $L^2(O_\rho)$ truncation error estimates:
\begin{equation*}
    \|\phi - \phi^{(L)}\|_{L^2(O_\rho)} \leq C_{\phi} \|\phi\|_{L^2(\partial O_d)}\left(\frac{\rho}{d}\right)^{\frac{L}{2}}, \quad 
    \|\psi - \psi^{(L)}\|_{L^2(O_\rho)} \leq C_{\psi} \|\psi\|_{L^2(\partial O_d)}\left(\frac{\rho}{d}\right)^{\frac{L}{2}},
\end{equation*}
where $C_{\phi}, C_{\psi} > 0$ are constants independent of $L$.

Consequently, under polar coordinates, the biharmonic function $u(r,\theta)$ can be expressed as:
\begin{equation*}
    u(r,\theta) = \sum_{n \in \mathbb{Z}} \xi_n\varphi_n + r^2\sum_{n \in \mathbb{Z}} \eta_n\varphi_n = \sum_{n \in \mathbb{Z}} \left( \xi_n\Phi_n^{(1)} + \eta_n\Phi_n^{(2)} \right).
\end{equation*}
Letting $\boldsymbol{\mu} = (\boldsymbol{\xi}, \boldsymbol{\eta})$ denote the combined coefficient vector, we have:
\begin{equation*}
    \|\boldsymbol{\mu} \|_{\ell^2}^2 = \|\boldsymbol{\xi} \|_{\ell^2}^2 + \|\boldsymbol{\eta} \|_{\ell^2}^2 \leq \frac{\rho}{d}\left( \|\phi\|_{L^2(\partial O_d)}^2 + \|\psi\|_{L^2(\partial O_d)}^2 \right).
\end{equation*}
% Taking the square root on both sides and utilizing the algebraic inequality $\sqrt{A^2 + B^2} \leq A + B$ for non-negative terms, we directly obtain the desired coefficient bound \eqref{eq:coefficient_bound}:
% \begin{equation*}
%     \|\boldsymbol{\mu}\|_{\ell^2} \leq \sqrt{\frac{\rho}{d}} \left( \|\phi\|_{L^2(\partial O_d)} + \|\psi\|_{L^2(\partial O_d)} \right).
% \end{equation*}

Next, we establish the truncation error estimate in the domain $O_\rho$. By the triangle inequality and noting that $r \leq \rho$ holds everywhere inside the disk $O_\rho$ (which implies $r^2 \leq \rho^2$), we have:
\begin{equation*}
\begin{aligned}
    \|u - u^{(L)}\|_{L^2(O_\rho)} 
    &\leq \|\phi - \phi^{(L)}\|_{L^2(O_\rho)} + \|r^2(\psi - \psi^{(L)})\|_{L^2(O_\rho)} \\
    &\leq \|\phi - \phi^{(L)}\|_{L^2(O_\rho)} + \rho^2 \|\psi - \psi^{(L)}\|_{L^2(O_\rho)}.
\end{aligned}
\end{equation*}
Substituting the individual harmonic truncation error bounds into the above inequality yields:
\begin{equation*}
\begin{aligned}
    \|u - u^{(L)}\|_{L^2(O_\rho)} 
    &\leq C_{\phi} \|\phi\|_{L^2(\partial O_d)}\left(\frac{\rho}{d}\right)^{\frac{L}{2}} + \rho^2 C_{\psi} \|\psi\|_{L^2(\partial O_d)}\left(\frac{\rho}{d}\right)^{\frac{L}{2}} \\
    &\leq C \left( \|\phi\|_{L^2(\partial O_d)} + \|\psi\|_{L^2(\partial O_d)} \right) \left(\frac{\rho}{d}\right)^{\frac{L}{2}},
\end{aligned}
\end{equation*}
where $C = \max(C_{\phi}, \rho^2 C_{\psi}) > 0$ is a constant independent of the truncation degree $L$. This completes the proof.
\end{proof}

\begin{lemma}
\label{lem:biharmonic-coefficient-to-interior}
For $K\in\mathbb N_0$, let
\begin{equation}
\label{eq:wK-definition}
    u_K(x)
    :=
    \sum_{n=-K}^{K}
    \left(
        \xi_n\varphi_n(x)
        +
        \eta_n|x|^2\varphi_n(x)
    \right),
\end{equation}
and define $
    \boldsymbol c
    :=
    \bigl(
        \xi_{-K},\ldots,\xi_K,
        \eta_{-K},\ldots,\eta_K
    \bigr)^T.$ 
Then
\begin{equation}
\label{eq:coefficient-to-L2}
    \|u_K\|_{L^2(O_\rho)}
    \leq
    C_\rho\|\boldsymbol c\|_{\ell^2},
\end{equation}
where $C_\rho>0$ is independent of $K$.
\end{lemma}

\begin{proof}
For $m=|n|$, direct integration in polar coordinates gives
\begin{equation}
    \|\varphi_n\|_{L^2(O_\rho)}^2
    =
    \int_0^\rho\int_0^{2\pi}
    \frac{r^{2m}}{2\pi\rho\,\rho^{2m}}
    r\,\mathrm{d}\theta\,\mathrm{d}r
    \nonumber
    =
    \frac{\rho}{2(m+1)},
\label{eq:phi-n-disk-norm}
\end{equation}
and
\begin{equation}
    \||x|^2\varphi_n\|_{L^2(O_\rho)}^2
    =
    \int_0^\rho\int_0^{2\pi}
    \frac{r^{2m+4}}{2\pi\rho\,\rho^{2m}}
    r\,\mathrm{d}\theta\,\mathrm{d}r
    \nonumber
    =
    \frac{\rho^5}{2(m+3)}.
\label{eq:r2-phi-n-disk-norm}
\end{equation}
Distinct Fourier modes are orthogonal in the angular variable.
Therefore,
\begin{equation}
    \left\|
        \sum_{n=-K}^{K}\xi_n\varphi_n
    \right\|_{L^2(O_\rho)}^2
    =
    \sum_{n=-K}^{K}
    \frac{\rho}{2(|n|+1)}|\xi_n|^2
    \nonumber
    \leq
    \frac{\rho}{2}
    \sum_{n=-K}^{K}|\xi_n|^2,
\label{eq:first-family-L2-bound}
\end{equation}
whereas
\begin{equation}
    \left\|
        \sum_{n=-K}^{K}\eta_n|x|^2\varphi_n
    \right\|_{L^2(O_\rho)}^2
    =
    \sum_{n=-K}^{K}
    \frac{\rho^5}{2(|n|+3)}|\eta_n|^2
    \nonumber
    \leq
    \frac{\rho^5}{6}
    \sum_{n=-K}^{K}|\eta_n|^2.
\label{eq:second-family-L2-bound}
\end{equation}
No orthogonality between the two different families is used.
By the triangle inequality,
\begin{equation*}
\begin{aligned}
    \|w_K\|_{L^2(O_\rho)}^2
    &\leq
    2
    \left\|
        \sum_{n=-K}^{K}\xi_n\varphi_n
    \right\|_{L^2(O_\rho)}^2
    +
    2
    \left\|
        \sum_{n=-K}^{K}\eta_n|x|^2\varphi_n
    \right\|_{L^2(O_\rho)}^2
    \nonumber\\
    &\leq
    \rho\sum_{n=-K}^{K}|\xi_n|^2
    +
    \frac{\rho^5}{3}
    \sum_{n=-K}^{K}|\eta_n|^2
    \nonumber
    \leq
    \max\left\{\rho,\frac{\rho^5}{3}\right\}
    \|\boldsymbol c\|_{\ell^2}^2.
\label{eq:wK-L2-square-bound}
\end{aligned}
\end{equation*}
This proves \eqref{eq:coefficient-to-L2}.
\end{proof}

The convergence theorem below follows directly from harmonic-mode truncation, the two-trace boundary-to-interior estimate, and the regularized collocation argument.

\begin{lemma}
\label{lem:biharmonic-boundary-stability}
Let $\Omega\subset\mathbb R^2$ be a bounded $C^4$ domain.
If $v\in H^2(\Omega)$ satisfies
\begin{equation*}
    \Delta^2v=0
    \qquad\text{in }\Omega,
\end{equation*}
then
\begin{equation*}
    \|v\|_{L^2(\Omega)}
    \leq
    C_{\mathrm{stab}}
    \left(
        \|v\|_{L^2(\partial\Omega)}
        +
        \|\partial_{n}v\|_{L^2(\partial\Omega)}
    \right),
\end{equation*}
where $C_{\mathrm{stab}}>0$ depends only on $\Omega$.
\end{lemma}

\begin{proof}
For arbitrary $F\in L^2(\Omega;\mathbb C)$, let $z$ be the solution
of the adjoint clamped problem
\begin{equation*}
\begin{cases}
    \Delta^2 z=F,
        &\text{in }\Omega,\\
    z=0,
        &\text{on }\partial\Omega,\\
    \partial_{n}z=0,
        &\text{on }\partial\Omega.
\end{cases}
\end{equation*}
By the regularity theory for the clamped biharmonic problem \cite[Corollary 2.21]{Gazzola_polyharmonic},
\begin{equation*}
    z\in H^4(\Omega;\mathbb C),
    \qquad
    \|z\|_{H^4(\Omega)}
    \leq
    C\|F\|_{L^2(\Omega)}.
\end{equation*}

We use the complex $L^2$-pairing
\begin{equation*}
    \langle f,g\rangle_{L^2(\Omega)}
    :=
    \int_\Omega f\,\overline g\,\mathrm{d}\boldsymbol x.
\end{equation*}
Since $\Delta^2v=0$ and
$z=\partial_{n}z=0$ on $\partial\Omega$, the
biharmonic Green identity gives
\begin{equation*}
\begin{aligned}
    \langle v,F\rangle_{L^2(\Omega)}
    =
    \int_\Omega v\,\overline{\Delta^2z}\,
        \mathrm{d}\boldsymbol x
    =
    \int_{\partial\Omega}
    \left[
        v\,\overline{\partial_{n}\Delta z}
        -
        (\partial_{n}v)\,
        \overline{\Delta z}
    \right]\,\mathrm{d}s .
\end{aligned}
\end{equation*}
Consequently,
\begin{equation*}
\begin{aligned}
    \left|
        \langle v,F\rangle_{L^2(\Omega)}
    \right|
    \leq{}
    \,
    \|v\|_{L^2(\partial\Omega)}
    \|\partial_{n}\Delta z\|_{L^2(\partial\Omega)}
    +
    \|\partial_{n}v\|_{L^2(\partial\Omega)}
    \|\Delta z\|_{L^2(\partial\Omega)}.
\end{aligned}
\end{equation*}
The trace theorem and the $H^4$-regularity estimate imply
\begin{equation*}
    \|\Delta z\|_{L^2(\partial\Omega)}
    +
    \|\partial_{n}\Delta z\|_{L^2(\partial\Omega)}
    \leq
    C\|F\|_{L^2(\Omega)}.
\end{equation*}
Therefore,
\begin{equation*}
    \left|
        \langle v,F\rangle_{L^2(\Omega)}
    \right|
    \leq
    C
    \left(
        \|v\|_{L^2(\partial\Omega)}
        +
        \|\partial_{n}v\|_{L^2(\partial\Omega)}
    \right)
    \|F\|_{L^2(\Omega)}.
\end{equation*}
Taking the supremum over all
$F\in L^2(\Omega;\mathbb C)$ satisfying
$\|F\|_{L^2(\Omega)}=1$ yields
\begin{equation*}
    \|v\|_{L^2(\Omega)}
    \leq
    C
    \left(
        \|v\|_{L^2(\partial\Omega)}
        +
        \|\partial_{n}v\|_{L^2(\partial\Omega)}
    \right),
\end{equation*}
which completes the proof.
\end{proof}

We now establish an error bound for the regularized numerical solution $u^{(c)}$ that accounts for basis truncation and boundary collocation.

\begin{theorem}
\label{thm:main_convergence}
Let $0<\rho<d$, let $O_\rho$ and $O_d$ be concentric disks
centered at the origin, and let $\Omega$ be a bounded $C^4$ domain
with $\overline{\Omega}\subset O_\rho$. Assume that the exact solution extends to a function $u\in H^2(O_d)$ that is biharmonic in $O_d$, with
Almansi decomposition
$u=\phi+r^2\psi$.
Let
\[
    u^{(c)}(\boldsymbol{x})
    =
    \boldsymbol{a}_*(\boldsymbol{x})\boldsymbol{h}
\]
be the regularized numerical solution obtained from a quasi-uniform
partition of $\partial\Omega$ into $N$ collocation arcs and
truncation degree $K$. Set
$ L=2K+1$
for the number of modes in each harmonic family and
$ M=2L=4K+2$
for the total number of basis functions, and assume that there are
constants $0<c_\alpha\leq C_\alpha$, independent of $K$ and $N$,
such that
\begin{equation*}
    c_\alpha N(\rho/d)^{2K}
    \leq
    \alpha
    \leq
    C_\alpha N(\rho/d)^{2K}.
\end{equation*}
Then
\begin{equation*}
    \|u-u^{(c)}\|_{L^2(\Omega)}
    \leq
    C
    \left(
        \|\phi\|_{L^2(\partial O_d)}
        +
        \|\psi\|_{L^2(\partial O_d)}
    \right)
    \left(
        \left(\frac{\rho}{d}\right)^{\frac{L}{2}}
        +
        \frac{1}{N}
    \right),
\end{equation*}
where $C>0$ is independent of $L$ and $N$ but may depend on
$\Omega$, $\rho$, $d$, and the quasi-uniformity and regularization
constants.
\end{theorem}

\begin{proof}
By definition in \eqref{eq:tikhonov_row}, the complex reconstruction can be written as
\begin{equation}
   \widetilde{u}^{(c)}(\boldsymbol{x})
    =
    \boldsymbol{a}_*(\boldsymbol{x})\boldsymbol{h}
    =
    \boldsymbol{b}(\boldsymbol{x})
    V^*(VV^*+\alpha I)^{-1}\boldsymbol{h},
\end{equation}
where $\boldsymbol{b}(\boldsymbol{x})\in \mathbb{C}^{1\times M}$ is the row vector of
biharmonic basis functions. Using the identity
$
    V^*(VV^*+\alpha I)^{-1}
    =
    (V^*V+\alpha I)^{-1}V^*,
 $
this expression becomes
\begin{equation}
   \widetilde{ u}^{(c)}(\boldsymbol{x})
    =
    \boldsymbol{b}(\boldsymbol{x})\boldsymbol{c}_*,
    \quad
    \mbox{with}\quad
    \boldsymbol{c}_*
    =
    (V^*V+\alpha I)^{-1}V^*\boldsymbol{h}.
\end{equation}
Since $u$ and $\boldsymbol{h}$ are real valued, $|u - \operatorname{Re}\widetilde z| \leq |u-z|$ for every $z \in \mathbb{C}$. Thus the residual and derivative estimates for the complex reconstruction also bound the reported approximation $u^{(c)}=\operatorname{Re}\widetilde{ u}^{(c)}$; for readability, we write $u^{(c)}$ in the remaining estimates.
The coefficient vector $\boldsymbol{c}_*\in\mathbb{C}^M$ is the
unique minimizer of the Tikhonov-regularized functional
\begin{equation}
    \boldsymbol{c}_*
    =
    \arg\min_{\boldsymbol{c}\in\mathbb{C}^M}
    \left\{
        \|\boldsymbol{h}-V\boldsymbol{c}\|^2
        +
        \alpha\|\boldsymbol{c}\|^2
    \right\}.
\end{equation}
Choosing the exact truncated coefficient vector
$\boldsymbol{\mu}=(\boldsymbol{\xi},\boldsymbol{\eta})$ from
Lemma~\ref{lem:biharmonic_polynomial_approximation} as a candidate
minimizer, the optimality of $\boldsymbol{c}_*$ yields
\begin{equation}
\label{eq:total_residuals_estimate}
    \|\boldsymbol{h}-V\boldsymbol{c}_*\|^2
    +
    \alpha\|\boldsymbol{c}_*\|^2
    \leq
    \|\boldsymbol{h}-V\boldsymbol{\mu}\|^2
    +
    \alpha\|\boldsymbol{\mu}\|^2.
\end{equation}
The term $\|\boldsymbol{h}-V\boldsymbol{\mu}\|^2$ represents the
interpolation error at the $N$ boundary collocation points:
\begin{equation}
\begin{aligned}
    \|\boldsymbol{h}-V\boldsymbol{\mu}\|^2
    ={}
    \sum_{j=1}^N
    |u(\boldsymbol{x}_j)-u^{(\mu)}(\boldsymbol{x}_j)|^2
    +
    \sum_{j=1}^N
    |\partial_{n}u(\boldsymbol{x}_j)
    -
    \partial_{n}u^{(\mu)}(\boldsymbol{x}_j)|^2.
\end{aligned}
\end{equation}
Let
$ e(\boldsymbol{x})
    =
    u(\boldsymbol{x})-u^{(\mu)}(\boldsymbol{x})$
be the biharmonic truncation error. Since both $u$ and
$u^{(\mu)}$ are biharmonic in $O_\rho$, we have
$\Delta^2e=0$. Choose a smooth domain $D$ with
$\Omega\Subset D\Subset O_\rho.$
The embedding $H^3(D)\hookrightarrow W^{1,\infty}(D)$ in two
dimensions and the standard interior estimate \cite[Theorem 2.20]{Gazzola_polyharmonic} for the biharmonic
operator give
\begin{equation*}
\begin{aligned}
    \|e\|_{L^\infty(\partial\Omega)}^2
    +&
    \|\partial_{n}e\|_{L^\infty(\partial\Omega)}^2
    \leq
    C_{\mathrm{emb}}\|e\|_{H^3(D)}^2
    \\
    &\leq
    C\|e\|_{L^2(O_\rho)}^2
    \leq
    C
    \left(
        \|\phi\|_{L^2(\partial O_d)}^2
        +
        \|\psi\|_{L^2(\partial O_d)}^2
    \right)
    \left(\frac{\rho}{d}\right)^{2K+1}.
\end{aligned}
\end{equation*}
where the last inequality is obtained by Lemma \ref{lem:biharmonic_polynomial_approximation}. Summing these pointwise errors over the $N$ collocation points yields
\begin{equation}
\label{eq:residuals_estimate}
\begin{aligned}
    \|\boldsymbol{h}-V\boldsymbol{\mu}\|^2
    \leq
    \sum_{j=1}^N
    \left(
        |e(\boldsymbol{x}_j)|^2
        +
        |\nabla e(\boldsymbol{x}_j)|^2
    \right)
    \leq
    CN
    \left(\frac{\rho}{d}\right)^{2K+1}
    \left(
        \|\phi\|_{L^2(\partial O_d)}^2
        +
        \|\psi\|_{L^2(\partial O_d)}^2
    \right).
\end{aligned}
\end{equation}

Substituting \eqref{eq:residuals_estimate} and the coefficient bound
\eqref{eq:coefficient_bound} into
\eqref{eq:total_residuals_estimate}, and using
$\alpha\asymp N(\rho/d)^{2K},$
we obtain
\begin{equation}
\label{eq:optimal_functional_bound}
\begin{aligned}
    \|\boldsymbol{h}-V\boldsymbol{c}_*\|^2
    +
    \alpha\|\boldsymbol{c}_*\|^2
    \leq
    CN
    \left(\frac{\rho}{d}\right)^{2K+1}
    \left(
        \|\phi\|_{L^2(\partial O_d)}^2
        +
        \|\psi\|_{L^2(\partial O_d)}^2
    \right).
\end{aligned}
\end{equation}
This bound implies two key inequalities:
\begin{equation}
\label{eq:resid_bound}
    \|\boldsymbol{h}-V\boldsymbol{c}_*\|^2
    \leq
    CN
    \left(\frac{\rho}{d}\right)^{2K+1}
    \left(
        \|\phi\|_{L^2(\partial O_d)}^2
        +
        \|\psi\|_{L^2(\partial O_d)}^2
    \right),
\end{equation}
and
\begin{equation}
\label{eq:coeff_bound_star}
\begin{aligned}
    \|\boldsymbol{c}_*\|^2
    \leq
    \frac{CN(\rho/d)^{2K+1}}{\alpha}
    \left(
        \|\phi\|_{L^2(\partial O_d)}^2
        +
        \|\psi\|_{L^2(\partial O_d)}^2
    \right) 
    \leq
    C
    \left(
        \|\phi\|_{L^2(\partial O_d)}^2
        +
        \|\psi\|_{L^2(\partial O_d)}^2
    \right).
\end{aligned}
\end{equation}

To estimate the continuous boundary error
$v(\boldsymbol{x})
    =
    u(\boldsymbol{x})-u^{(c)}(\boldsymbol{x}),$
we partition the boundary as
$\partial\Omega
    =
    \bigcup_{k=1}^N\Gamma_k$,
such that $\boldsymbol{x}_k\in\Gamma_k$ and
$|\Gamma_k|\leq C/N$. On each segment $\Gamma_k$, the triangle
inequality yields
\begin{equation}
\label{eq:error_bound_in_small_domain}
\begin{aligned}
    \int_{\Gamma_k}|v(\boldsymbol{x})|^2
        \,\mathrm{d}\boldsymbol{x}
    \leq{}
    2\int_{\Gamma_k}
    |v(\boldsymbol{x})-v(\boldsymbol{x}_k)|^2
        \,\mathrm{d}\boldsymbol{x}
        +
    2\int_{\Gamma_k}|v(\boldsymbol{x}_k)|^2
        \,\mathrm{d}\boldsymbol{x}.
\end{aligned}
\end{equation}
Summing the second term in
\eqref{eq:error_bound_in_small_domain} over all segments and applying
\eqref{eq:resid_bound} gives
\begin{equation}
\label{eq:discrete_part}
\begin{aligned}
    \sum_{k=1}^N
    \int_{\Gamma_k}|v(\boldsymbol{x}_k)|^2
        \,\mathrm{d}\boldsymbol{x}
    \leq
    \frac{C}{N}
    \|\boldsymbol{h}-V\boldsymbol{c}_*\|^2
    \leq
    C
    \left(\frac{\rho}{d}\right)^{2K+1}
    \left(
        \|\phi\|_{L^2(\partial O_d)}^2
        +
        \|\psi\|_{L^2(\partial O_d)}^2
    \right).
\end{aligned}
\end{equation}
For the first term in
\eqref{eq:error_bound_in_small_domain}, we split it via the
oscillations:
\begin{equation*}
\begin{aligned}
    \int_{\Gamma_k}
    |v(\boldsymbol{x})-v(\boldsymbol{x}_k)|^2
        \,\mathrm{d}\boldsymbol{x}
    \leq{}
    2\int_{\Gamma_k}
    |u(\boldsymbol{x})-u(\boldsymbol{x}_k)|^2
        \,\mathrm{d}\boldsymbol{x}
    +
    2\int_{\Gamma_k}
    |u^{(c)}(\boldsymbol{x})
    -
    u^{(c)}(\boldsymbol{x}_k)|^2
        \,\mathrm{d}\boldsymbol{x}.
\end{aligned}
\end{equation*}
Using the mean-value theorem along an arc-length parametrization,
with
$ s(\boldsymbol{x})
    \leq
    |\Gamma_k|
    \leq
    C/N,$
we have
\begin{equation*}
\begin{aligned}
    \int_{\Gamma_k}
    |u(\boldsymbol{x})-u(\boldsymbol{x}_k)|^2
        \,\mathrm{d}\boldsymbol{x}
    \leq
    |\Gamma_k|^3
    \max_{\boldsymbol{x}\in\Gamma_k}
    |\nabla u(\boldsymbol{x})|^2
    \leq
    \frac{C}{N^3}\|u\|_{L^2(O_d)}^2
    \leq
    \frac{C}{N^3}
    \left(
        \|\phi\|_{L^2(\partial O_d)}^2
        +
        \|\psi\|_{L^2(\partial O_d)}^2
    \right),
\end{aligned}
\end{equation*}
where the second inequality follows from the interior derivative
estimate on $\Omega\Subset O_d$. The last inequality follows by writing $u = \phi + r^2\psi$ and applying Lemma \ref{lem: interior_stability_estimate} to the harmonic components
on $O_d$.
Similarly, an interior estimate for $u^{(c)}$ gives
\begin{equation*}
    \max_{\boldsymbol{x}\in\partial\Omega}
    |\nabla u^{(c)}(\boldsymbol{x})|^2
    \leq
    C\|u^{(c)}\|_{L^2(O_\rho)}^2.
\end{equation*}
The uniform coefficient bound \eqref{eq:coeff_bound_star} and
Lemma~\ref{lem:biharmonic-coefficient-to-interior} imply
\begin{equation*}
\begin{aligned}
    \|u^{(c)}\|_{L^2(O_\rho)}^2
    \leq
    C_\rho^2\|\boldsymbol{c}_*\|^2
    \leq
    C
    \left(
        \|\phi\|_{L^2(\partial O_d)}^2
        +
        \|\psi\|_{L^2(\partial O_d)}^2
    \right).
\end{aligned}
\end{equation*}
Consequently,
\begin{equation*}
\begin{aligned}
    \int_{\Gamma_k}
    |u^{(c)}(\boldsymbol{x})
    -
    u^{(c)}(\boldsymbol{x}_k)|^2
        \,\mathrm{d}\boldsymbol{x}
    \leq
    \frac{C}{N^3}
    \|u^{(c)}\|_{L^2(O_\rho)}^2
    \leq
    \frac{C}{N^3}
    \left(
        \|\phi\|_{L^2(\partial O_d)}^2
        +
        \|\psi\|_{L^2(\partial O_d)}^2
    \right).
\end{aligned}
\end{equation*}
Summing these components over all $k=1,\ldots,N$ yields
\begin{equation*}
\begin{aligned}
    \|u-u^{(c)}\|_{L^2(\partial\Omega)}^2
    \leq
    C
    \left(
        \left(\frac{\rho}{d}\right)^{2K+1}
        +
        \frac{1}{N^2}
    \right)
    \left(
        \|\phi\|_{L^2(\partial O_d)}^2
        +
        \|\psi\|_{L^2(\partial O_d)}^2
    \right).
\end{aligned}
\end{equation*}

For the Neumann trace, let $\gamma$ be a unit-speed parametrization
of an arc, let $\tau=\gamma'$, and let $n$ be the unit normal. For
$v=u-u^{(c)}$ and
$w(s)=\partial_n v(\gamma(s)),$ 
differentiation along the boundary gives
\begin{equation*}
    w'(s)
    =
    \bigl(D^2v(\gamma(s))\,\tau(s)\bigr)\cdot n(s)
    +
    \nabla v(\gamma(s))\cdot n'(s).
\end{equation*}
Because $\partial\Omega$ is $C^4$, its curvature is bounded, and
hence
\begin{equation*}
    |w'(s)|
    \leq
    C_\Omega
    \left(
        |D^2v(\gamma(s))|
        +
        |\nabla v(\gamma(s))|
    \right).
\end{equation*}
Interior derivative estimates on $\Omega\Subset O_\rho$, the
harmonic stability estimate for $u$, and
\eqref{eq:coeff_bound_star}, together with
Lemma \ref{lem:biharmonic-coefficient-to-interior} for $u^{(c)}$,
imply
\begin{equation*}
    \|w'\|_{L^\infty(\partial\Omega)}^2
    \leq
    C
    \left(
        \|\phi\|_{L^2(\partial O_d)}^2
        +
        \|\psi\|_{L^2(\partial O_d)}^2
    \right).
\end{equation*}
Thus the same arc-length argument used for the Dirichlet trace
bounds the sum of the within-arc oscillations by $CN^{-2}$ times the right-hand side above. The Neumann components of
$\boldsymbol{h}-V\boldsymbol{c}_*$ control the collocation values of
$w$, so \eqref{eq:resid_bound} yields
\begin{equation*}
\begin{aligned}
    \|\partial_n(u-u^{(c)})\|_{L^2(\partial\Omega)}^2
    \leq
    C
    \left(
        \left(\frac{\rho}{d}\right)^{2K+1}
        +
        \frac{1}{N^2}
    \right)
    \left(
        \|\phi\|_{L^2(\partial O_d)}^2
        +
        \|\psi\|_{L^2(\partial O_d)}^2
    \right).
\end{aligned}
\end{equation*}

Finally, invoking the biharmonic stability estimate in
Lemma~\ref{lem:biharmonic-boundary-stability}, we obtain
\begin{equation*}
\begin{aligned}
    \|u-u^{(c)}\|_{L^2(\Omega)}
    \leq
    C_{\mathrm{stab}}
    \bigl(
        \|u-u^{(c)}\|_{L^2(\partial\Omega)}
        +
        \|\partial_n(u-u^{(c)})\|_{L^2(\partial\Omega)}
    \bigr).
\end{aligned}
\end{equation*}
Taking square roots and using $L=2K+1$ completes the proof.
\end{proof}

The component-dependent quantity appearing in Theorem \ref{thm:main_convergence}
can be controlled directly by the Sobolev norm of the extended exact solution.
Combining the Sobolev-stable Almansi decomposition in Theorem \ref{thm:stable-almansi-decomposition} yields the following consequence.

\begin{corollary}[Error estimate in terms of the extended solution]
\label{cor:main-convergence-solution-norm}
Under the assumptions of Theorem~\ref{thm:main_convergence}, the regularized
numerical solution satisfies
\begin{equation}
\label{eq:main-convergence-solution-norm}
    \|u-u^{(c)}\|_{L^2(\Omega)}
    \leq
    C
    \|u\|_{H^2(O_d)}
    \left[
        \left(\frac{\rho}{d}\right)^{L/2}
        +
        \frac{1}{N}
    \right],
\end{equation}
where \(C>0\) is independent of \(u\), \(L\), and \(N\), but may depend on
\(\Omega\), \(\rho\), \(d\), and the quasi-uniformity and regularization
constants.
\end{corollary}

\section{Numerical Experiments}
\label{sec:experiments}

 \providecommand{\figureplaceholder}[1]{%
   \fbox{\begin{minipage}[c][0.82in][c]{0.91\textwidth}
   \centering\small #1
   \end{minipage}}}
 \providecommand{\revisioncolor}{}

In this section, we present numerical experiments in complex regions to demonstrate the effectiveness of the proposed LbNM. 

\subsection{Experimental Setup}
\label{sec:setup}

We first describe the geometries, discretization parameters, and computing environment.

We consider one rounded-kite domain and one flower-shaped domain. Their parametric boundaries are defined in Table \ref{tab:different-domains}. The expansion centers are $(0.1,0)$ for the rounded kite in both test cases and $(0,0)$ for the flower.
\begin{table}[htbp]
\revisioncolor
		\centering
	\caption{Parametric definitions of the computational domains.}
		\label{tab:different-domains}
	\begin{tabular}{cc}
		\toprule
		Shape & Parameter equation \\
		\midrule
		Rounded kite   & 
		\begin{tabular}[t]{@{}c@{}}
			$x = 0.5 \cos t + 0.3 \cos(2t) - 0.2, \quad t \in [0, 2\pi]$, \\
			$y = 0.6 \sin t, \quad t \in [0, 2\pi]$
		\end{tabular} \\
		\midrule
		Flower-shape & 
		\begin{tabular}[t]{@{}c@{}}
			$x = 0.5 \cos t - 0.1 \cos 6t \cos t, \quad t \in [0, 2\pi]$, \\
			$y = 0.5 \sin t - 0.1 \cos 6t \sin t, \quad t \in [0, 2\pi]$
		\end{tabular} \\
		\bottomrule
	\end{tabular}
\end{table}

In our implementation, $N$ boundary collocation points $\{x_j\}_{j=1}^{N}$ are uniformly distributed with respect to the parameter $t \in [0, 2\pi)$. Since the biharmonic boundary value problem involves two distinct boundary conditions (i.e., Dirichlet and Neumann data), the collocation matrix $V$ has $2N$ rows.

Let $K$ denote the maximum truncation degree of the harmonic polynomials. Each harmonic family contains $L=2K+1$ modes, and the two-family biharmonic hypothesis space therefore has $M=2L=4K+2$ columns. The trace matrix is factorized by an economy QR decomposition followed by an SVD of the triangular factor, and the coefficients are computed using the Tikhonov-regularized coefficients filter $\sigma_j/(\sigma_j^2+\alpha)$. 

Throughout this section, $E_\Omega$ denotes the relative discrete $\ell^2$ error on the reported interior evaluation set. On an independent boundary validation grid $X_{\partial\Omega}^{\mathrm{val}}$, the Dirichlet and normal-derivative trace errors are reported separately as
\[
E_u^{\mathrm{val}}
=\frac{\bigl(\sum_{x\in X_{\partial\Omega}^{\mathrm{val}}}|u_{\mathrm{LbNM}}(x)-u(x)|^2\bigr)^{1/2}}
{\bigl(\sum_{x\in X_{\partial\Omega}^{\mathrm{val}}}|u(x)|^2\bigr)^{1/2}}.
\]
\[
E_{\partial_nu}^{\mathrm{val}}
=\frac{\bigl(\sum_{x\in X_{\partial\Omega}^{\mathrm{val}}}|\partial_{n}u_{\mathrm{LbNM}}(x)-\partial_{n}u(x)|^2\bigr)^{1/2}}
{\bigl(\sum_{x\in X_{\partial\Omega}^{\mathrm{val}}}|\partial_{n}u(x)|^2\bigr)^{1/2}}.
\]

\subsubsection{Computing Environment}

All algorithms were implemented in MATLAB R2021a and run on a macOS computer with an Apple M4 processor and 16 GB of memory. Reported CPU times are wall-clock measurements from this platform.

\subsection{Test Case 1: Smooth Solution}
\label{subsec:smooth_benchmark}

We first consider a benchmark problem with a smooth analytical solution. The exact solution is prescribed as
\begin{equation}
\label{eq:entire_biharmonic_benchmark}
	u(x,y) = e^x\cos y + (x^2 + y^2)e^y\cos x,
\end{equation}
which satisfies $\Delta^2 u = 0$ on $\mathbb R^2$. In particular, it admits a biharmonic continuation to every concentric disk $O_d$ and thus meets the continuation hypothesis of Theorem \ref{thm:main_convergence}. 

\begin{table}[htbp]
	\revisioncolor
		\centering
		\caption{Discretization parameters for the computational domains.}
	\label{tab:parameters}
		\begin{tabular}{lccccc}
			\toprule
			Domain Shape & $N$ & $K$ & $M=4K+2$ & $\rho$ & Offine build (s)\\
			\midrule
			Rounded kite & 300 & 125 & 502 & 0.91 & 0.49 \\
			Flower       & 388 & 100 & 402 & 0.61 & 0.47 \\
			\bottomrule
		\end{tabular}
\end{table}

For the theorem-compatible parameter schedule we use $d=1.456$ for the rounded kite and $d=1.25$ for the flower. The numerical experiments use the discretization parameters specified in Table \ref{tab:parameters}.

\begin{figure}[htbp]
\centering
\includegraphics[width=\textwidth]
{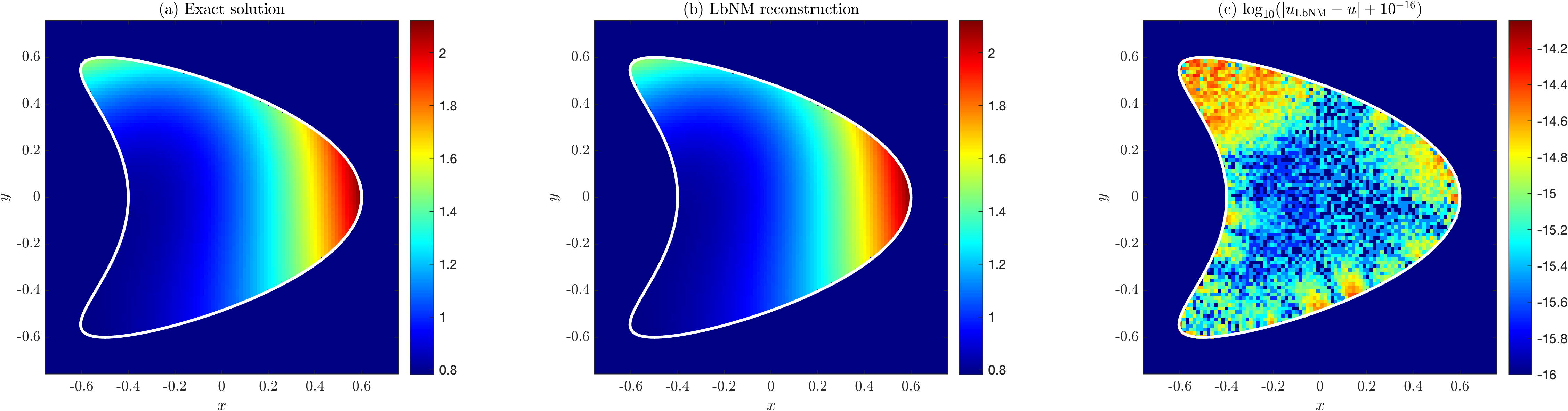}
\caption{Numerical performance of the proposed LbNM for the smooth solution on the rounded kite-shaped domain: (a) analytical exact solution, (b) LbNM prediction, and (c) spatial distribution of the pointwise absolute error on the $\log_{10}$ scale. The solid white curve represents the physical boundary.}
\label{fig:kite-entire-domain}
\end{figure}

\begin{figure}[htbp]
\centering
\includegraphics[width=\textwidth]
{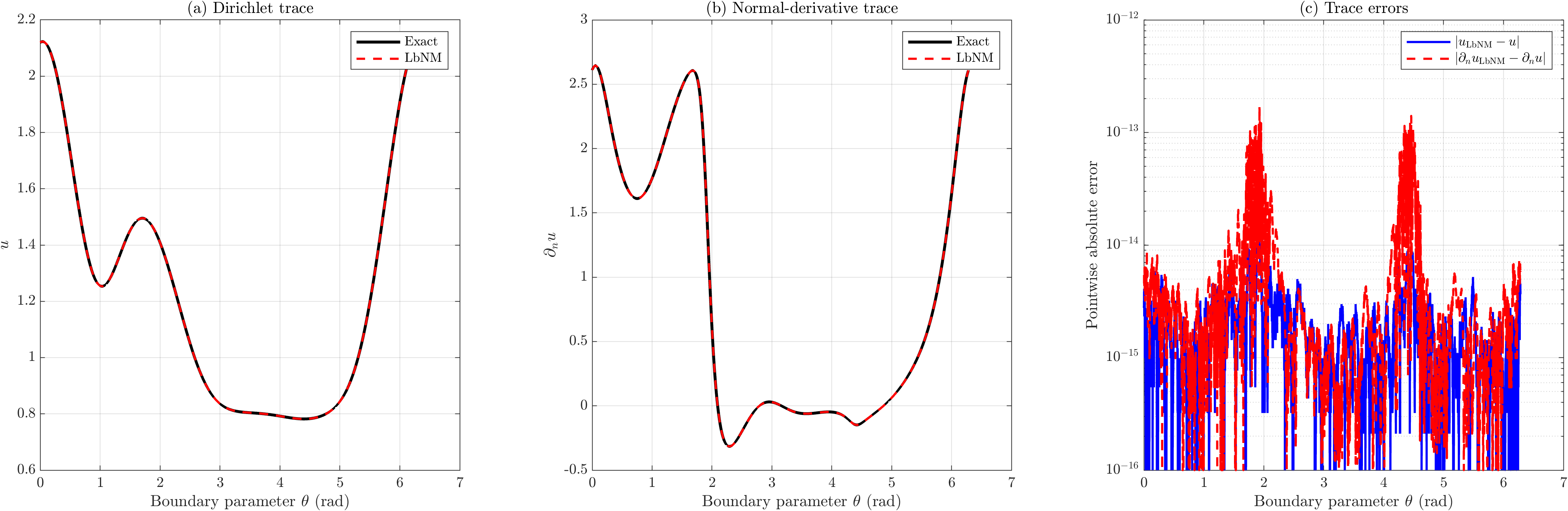}
\caption{Boundary numerical results for the smooth solution on the rounded kite-shaped domain: (a) boundary profile of $u$, (b) boundary profile of the normal derivative $\partial_{n}u$, and (c) pointwise absolute boundary errors for both $u$ and $\partial_{n}u$ plotted against the boundary angle $\theta$.}
\label{fig:kite-entire-boundary}
\end{figure}

\begin{figure}[htbp]
\centering
\includegraphics[width=\textwidth]
{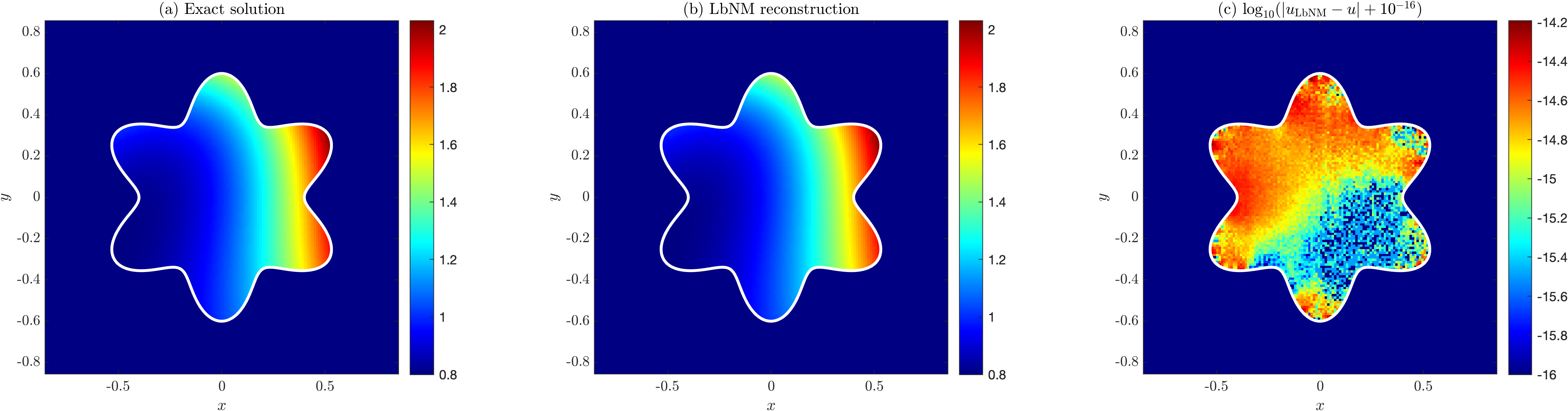}
\caption{Numerical performance of the proposed LbNM for the smooth solution on the flower-shaped domain: (a) analytical exact solution, (b) LbNM prediction, and (c) spatial distribution of the pointwise absolute error on the $\log_{10}$ scale. The solid white curve represents the physical boundary.}
\label{fig:flower-entire-domain}
\end{figure}

\begin{figure}[H]
\centering
\includegraphics[width=\textwidth]
{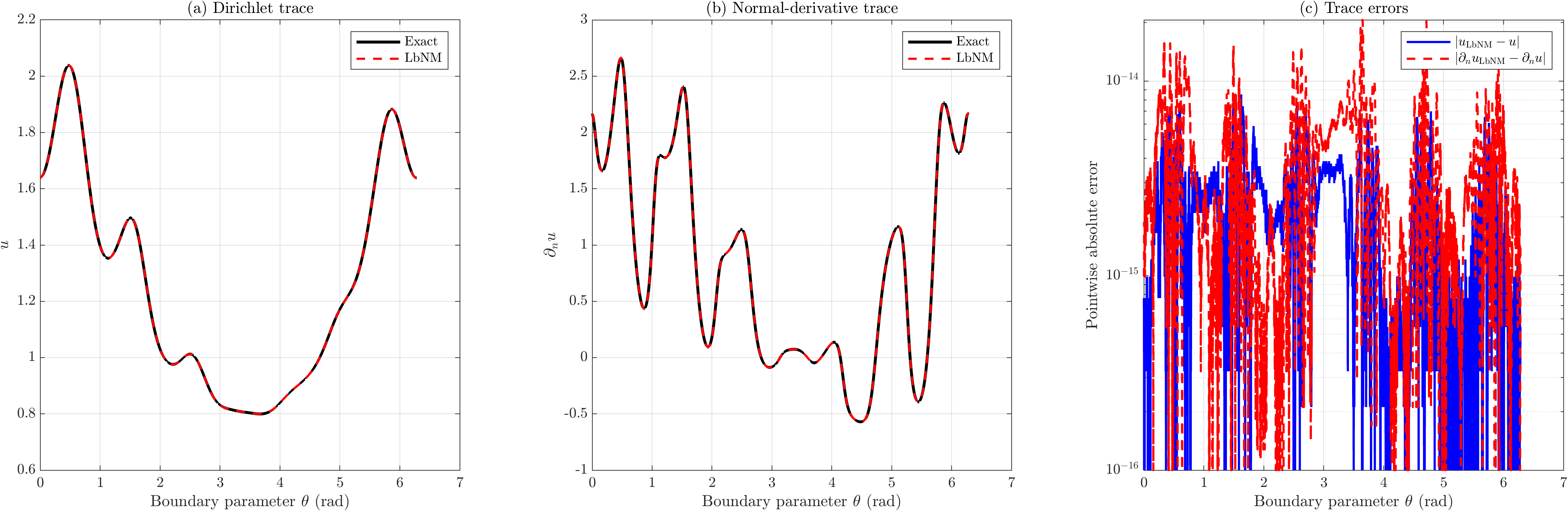}
\caption{Boundary numerical results for the smooth solution on the flower-shaped domain: (a) boundary profile of $u$, (b) boundary profile of the normal derivative $\partial_{n}u$, and (c) pointwise absolute boundary errors for both $u$ and $\partial_{n}u$ plotted against the boundary angle $\theta$.}
\label{fig:flower-entire-boundary}
\end{figure}

Figures \ref{fig:kite-entire-domain}--\ref{fig:flower-entire-boundary} compare the exact and reconstructed solutions and display pointwise and independently evaluated boundary errors.

\begin{table}[htbp]
	\revisioncolor
		\centering
		\caption{Relative numerical errors and CPU times for the smooth solution.}
	\label{tab:smooth_results}
	\small
	\begin{tabular}{lcccc}
		\toprule
		Domain Shape & $E_\Omega$ & $E_u^{\mathrm{val}}$ & $E_{\partial_nu}^{\mathrm{val}}$ & Online time (s) \\
		\midrule
		Rounded kite & $9.0446\times10^{-16}$ & $1.9700\times10^{-15}$
& $1.3377\times10^{-14}$
& $7.52\times10^{-4}$ \\
		Flower       & $1.4264\times10^{-15}$ & $1.9204\times10^{-15}$
& $4.0886\times10^{-15}$
& $8.05\times10^{-4}$ \\
		\bottomrule
	\end{tabular}
\end{table}

Table \ref{tab:smooth_results} shows errors at the double-precision floor. Comparable accuracy is obtained on both nonconvex geometries without evaluating either boundary error on the collocation grid.

\subsection{Test Case 2: Near-Boundary Singularity}
\label{subsec:singular_benchmark}
In this example, we verify that the operator learned in Test 1 applies to the near-boundary singularity solutions. For example, 
we consider a benchmark based on the two-dimensional biharmonic fundamental solution,
\begin{equation}
	u(\boldsymbol{x}) = |\boldsymbol{x} - \boldsymbol{x}_0|^2 \ln |\boldsymbol{x} - \boldsymbol{x}_0|,
\end{equation}
where $\boldsymbol{x}_0\in\mathbb{R}^2\setminus\overline\Omega$ is the source singularity. Specifically, $\boldsymbol{x}_0=(-0.60,0.00)$ for the rounded kite and $\boldsymbol{x}_0=(0.59,0.00)$ for the flower. 
%Because each source lies outside the physical domain, the exact solution is smooth and biharmonic in a neighborhood of $\overline\Omega$. This local fact is sufficient to define a regular boundary value problem, but it is weaker than the global concentric-disk continuation hypothesis of Theorem \ref{thm:main_convergence}. For the rounded kite, $|\boldsymbol x_0-\boldsymbol c|=0.70<\rho=0.91$ with $\boldsymbol c=(0.1,0)$; for the flower, $|\boldsymbol x_0|=0.59<\rho=0.61$. Thus both singularities lie inside their corresponding fictitious disks, and no admissible concentric extension radius $d>\rho$ can avoid either source. 

 Table \ref{tab:singular_results} reports the domain error together with the independently evaluated Dirichlet and normal-derivative errors. These results show that the reconstruction can remain highly accurate. 

\begin{table}[htbp]
	\revisioncolor
		\centering
		\caption{Relative numerical errors and CPU times for the singular solution.}
	\label{tab:singular_results}
	\small
	\begin{tabular}{lcccc}
		\toprule
		Domain Shape & $E_\Omega$ & $E_u^{\mathrm{val}}$ & $E_{\partial_nu}^{\mathrm{val}}$ & Online time (s) \\
		\midrule
		Rounded kite & $4.1839\times10^{-10}$ & $1.6100\times10^{-9}$
& $1.3752\times10^{-8}$
& $7.33\times10^{-4}$ \\
		Flower       & $7.7781\times10^{-11}$ & $2.0718\times10^{-10}$
& $3.0106\times10^{-11}$
& $8.39\times10^{-4}$ \\
		\bottomrule
	\end{tabular}
\end{table}

\begin{figure}[htbp]
	\revisioncolor
		\centering
		\includegraphics[width=\textwidth]{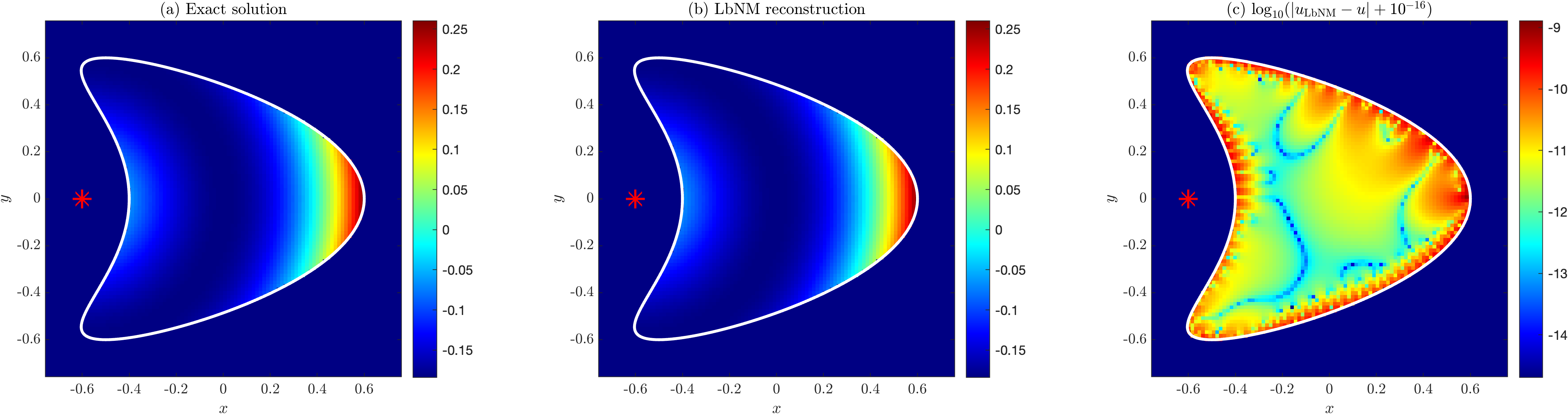}
	\caption{Numerical results of the LbNM  on a kite-shaped domain: (a) exact solution, (b) LbNM prediction, and (c) $\log_{10}$ domain error. The white solid line represents the domain boundary, and the red asterisk ($\ast$) denotes the location of the external source point.}
	\label{fig:LbNM_Kite_Fundamental_DomainResults}
\end{figure}

\begin{figure}[htbp]
	\revisioncolor
		\centering
		\includegraphics[width=\textwidth]{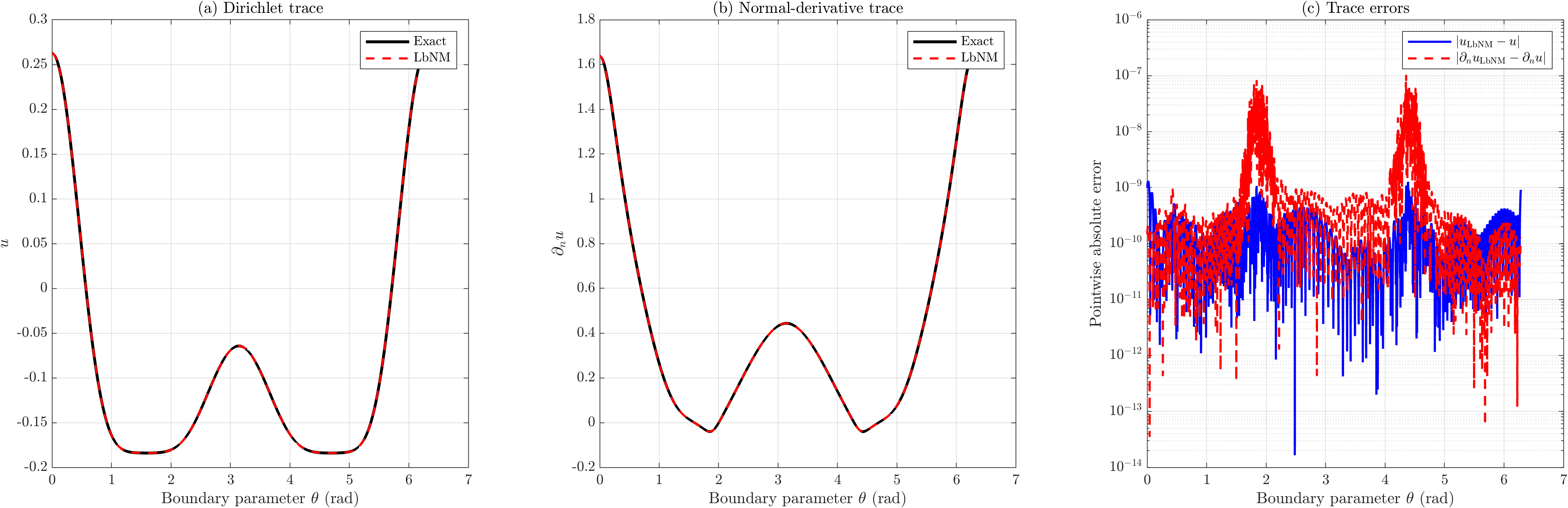}
	\caption{Boundary numerical results of the LbNM on the rounded kite-shaped domain: (a) boundary profile of $u$, (b) normal derivative $\partial_{n}u$ on the boundary, and (c) absolute boundary errors for both $u$ and $\partial_{n}u$. All results are plotted as functions of the boundary parameter $\theta$ (in radians).}
	\label{fig:LbNM_Kite_Fundamental_BoundaryResults}
\end{figure}

\begin{figure}[htbp]
	\revisioncolor
		\centering
		\includegraphics[width=\textwidth]{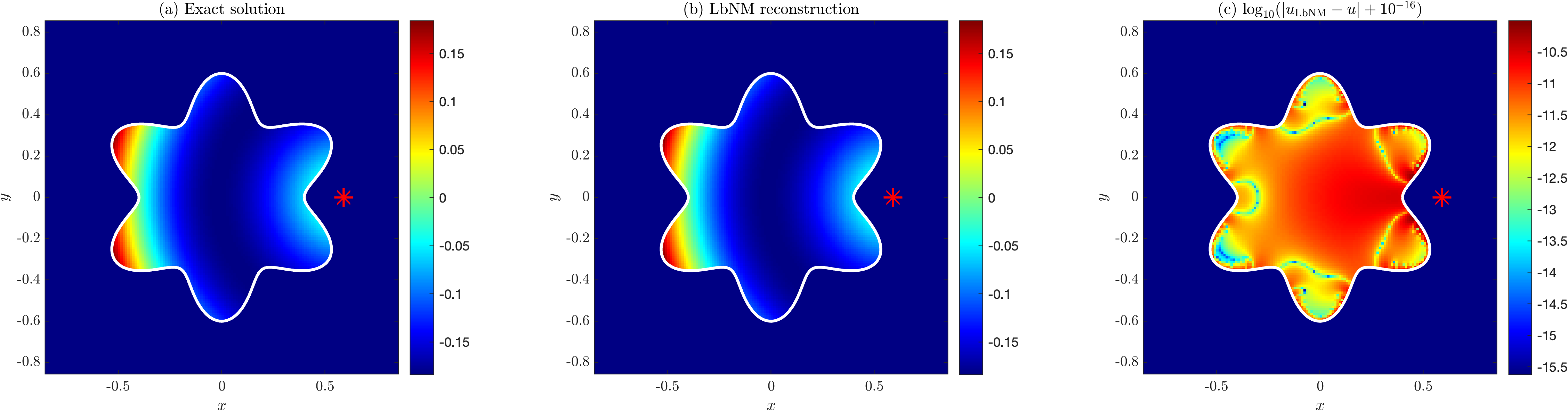}
	\caption{Numerical results of the LbNM on a flower-shaped domain: (a) exact solution, (b) LbNM prediction, and (c) $\log_{10}$ domain error. The white solid line represents the domain boundary, and the red asterisk ($\ast$) denotes the location of the external source point.}
	\label{fig:LbNM_Flower_Fundamental_DomainResults}
\end{figure}

\begin{figure}[htbp]
	\revisioncolor
		\centering
		\includegraphics[width=\textwidth]{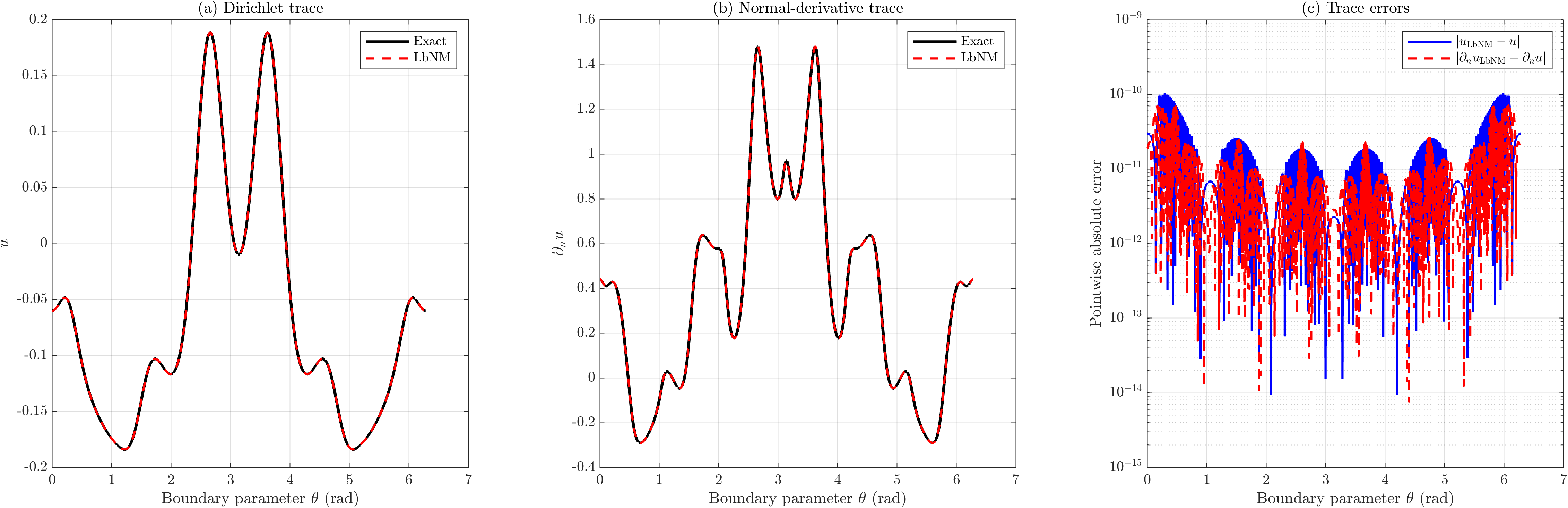}
	\caption{Boundary numerical results of the LbNM on the flower-shaped domain: (a) boundary profile of $u$, (b) normal derivative $\partial_{n}u$ on the boundary, and (c) absolute boundary errors for both $u$ and $\partial_{n}u$. All results are plotted as functions of the boundary parameter $\theta$ (in radians).}
	\label{fig:LbNM_Flower_Fundamental_BoundaryResults}
\end{figure}

Figures \ref{fig:LbNM_Kite_Fundamental_DomainResults}--\ref{fig:LbNM_Flower_Fundamental_BoundaryResults} show that the reconstructed fields capture the spatial variation induced by the nearby sources without visible spurious oscillations.

\subsection{Offline Study of the Approximation-Space Dimension}
\label{subsec:k_resolution}

\begin{figure}[htbp]
\centering
\includegraphics[width=0.5\textwidth]
{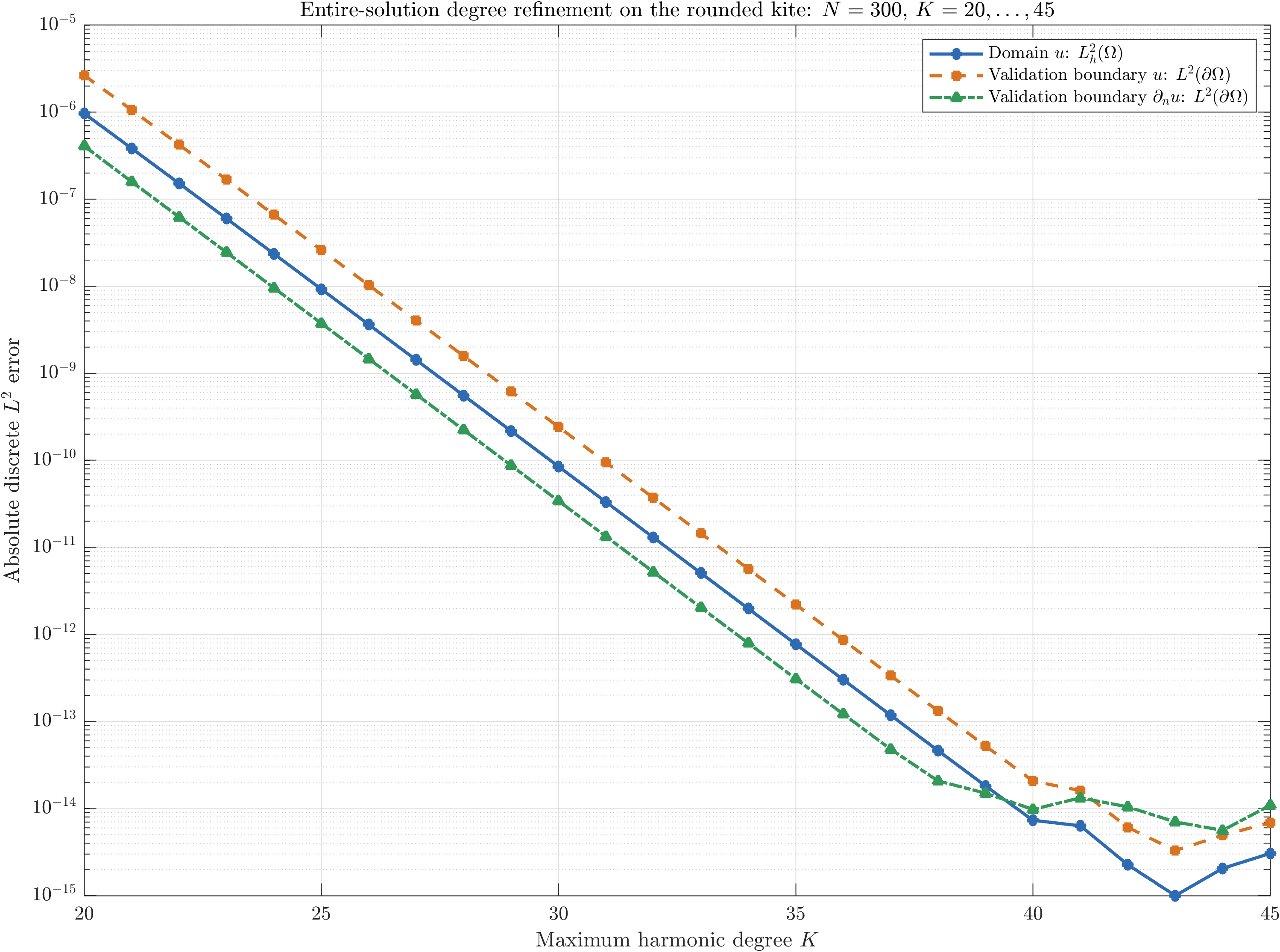}
\caption{Rounded-kite degree refinement for the entire solution with fixed $N=300$. }
\label{fig:kite-entire-k-refinement}
\end{figure}

\begin{figure}[htbp]
\centering
\includegraphics[width=0.5\textwidth]
{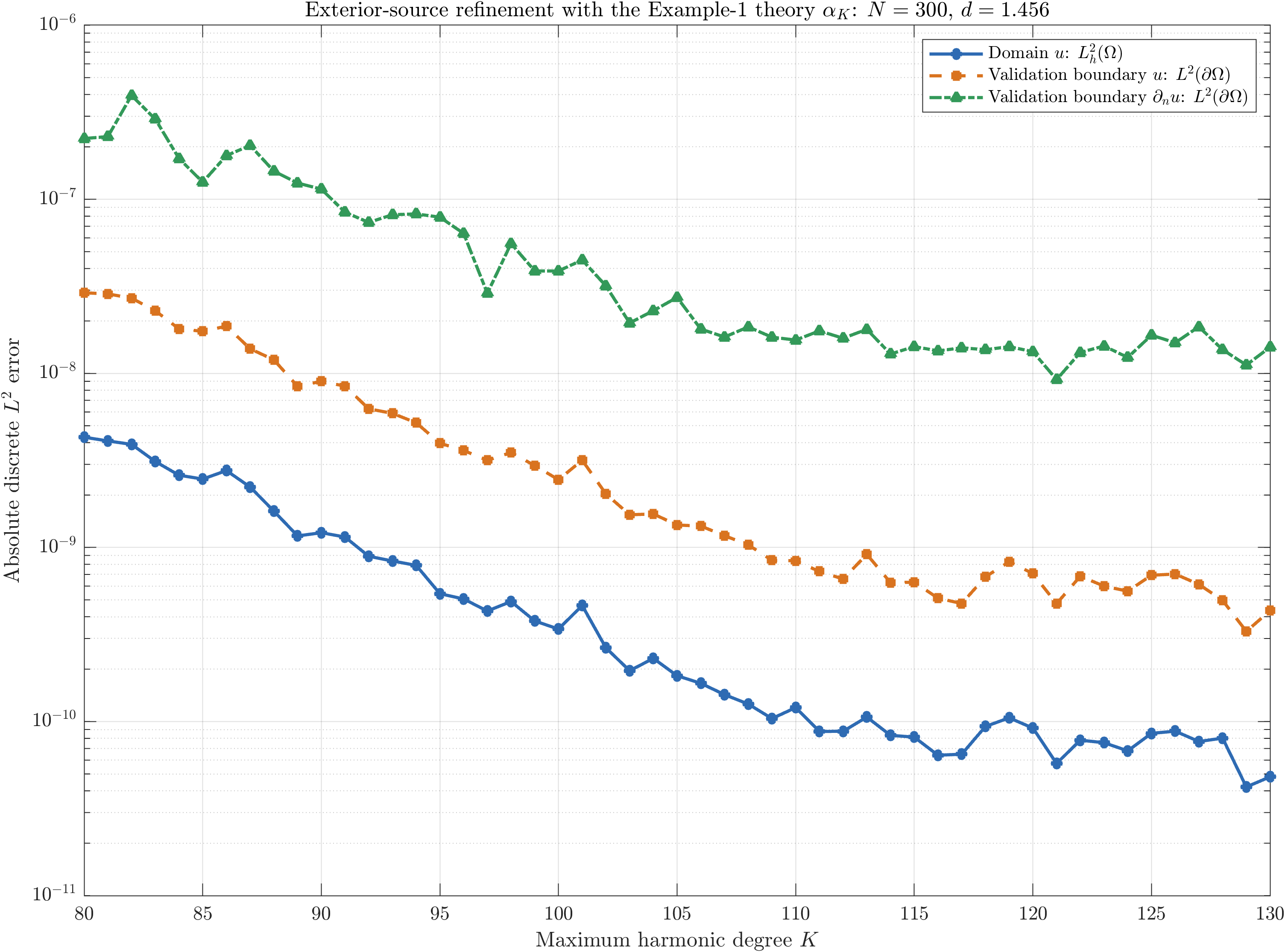}
\caption{Rounded-kite degree refinement for the near-boundary
exterior-source solution with fixed $N=300$. }
\label{fig:kite-exterior-k-refinement}
\end{figure}

The preceding two test cases use the same fixed production operator for each geometry, with Test Case 2 changing only the boundary-data vector. We additionally perform an offline $K$-refinement study on the rounded-kite domain to examine the resolution requirements of the two solutions. Each value of $K$ defines a newly constructed auxiliary operator and is therefore independent of the online reuse experiment. The results in Figure \ref{fig:kite-entire-k-refinement} and \ref{fig:kite-exterior-k-refinement} show that the entire solution is resolved with $K=20,\ldots,45$, whereas the near-boundary singularity requires the higher range $K=80,\ldots,130$ to achieve a comparable error reduction.

\section{Conclusion}
\label{sec:conclusions}

We have extended the Learning-Based Numerical Method to two-dimensional boundary value problems governed by the homogeneous biharmonic equation. The  method reconstructs a finite-dimensional boundary-to-solution operator from  samples generated by the Almansi decomposition and enforces both the solution and normal-derivative traces. Because the sample functions are biharmonic, the reconstructed approximation satisfies the governing equation by construction. The analysis provides stability and convergence estimates, while the numerical tests show high accuracy on two nonconvex geometries for both smooth data and data generated by nearby exterior singularities

The numerical results for solutions induced by singularities located near the boundary are particularly encouraging. Although such solutions exhibit increasingly localized variations as the singularity approaches the domain, the proposed method remains capable of producing accurate approximations. This observation motivates the development of adaptive sample-enrichment strategies. In future work, localized or singularity-informed sample solutions may be incorporated incrementally into the existing sample space, allowing global smooth modes and local singular features to be resolved within a unified framework.

\section*{Acknowledgments}
This work was supported by the China Postdoctoral Science Foundation [ Grant Number 2026M793352].

\section*{Declaration of generative AI and AI-assisted technologies in the manuscript preparation process.}
Statement: During the preparation of this work, the authors used Codex (OpenAI) to engage in interactive discussions about the manuscript,  and improve the clarity of the English text. After using this tool, the authors reviewed and edited the content as needed and take full responsibility for the content of the published article.

%This work of M. Bai is support in part by China Postdoctoral Science Foundation under Grant (2026M793352).

\bibliographystyle{plain}
\bibliography{references}

\end{document}